\documentclass[12pt,reqno]{amsart}
\usepackage{amssymb,latexsym,amsmath}
\usepackage{amsthm, color}
\usepackage[all]{xy}

\usepackage[pagebackref,colorlinks,
linkcolor=blue,
citecolor=blue,urlcolor=black,
hypertexnames=false]{hyperref}
\newcommand{\minimatrix}[4]
{\bigl(\begin{smallmatrix} 
#1 & #2 \\ #3 & #4\end{smallmatrix} \bigr)}
\newcommand{\cB}{\mathcal{B}}
\newcommand{\Cf}{\mathrm{C}}
\newcommand{\Co}{\mathcal{C}}

\newcommand{\cF}{\mathcal{F}}
\newcommand{\Hi}{\mathcal{H}}

\newcommand{\cK}{\mathcal{K}}
\newcommand{\cL}{\mathcal{L}}
\newcommand{\cM}{{\mathcal{M}}}
\newcommand{\On}{{\mathcal{O}_n}}

\newcommand{\Oinf}{{\mathcal{O}_\infty}}
\newcommand{\ot}{\otimes}

\newcommand{\cS}{\mathcal{S}}

\newcommand{\cY}{\mathcal{Y}}
\newcommand{\C}{\mathbb{C}}
\newcommand{\R}{\mathbb{R}}
\newcommand{\Z}{\mathbb{Z}}
\newcommand{\N}{\mathbb{N}}
\newcommand{\pf}
{\noindent{\mbox{\textbf{Proof}.\,}} }
\newcommand{\cf}{\textrm{cf.~}}
\newcommand{\eg}{\textrm{e.g.~}}
\newcommand{\ie}{\textrm{i.e., }}
\newcommand{\cst}{\textit{C}*}

\newcommand{\diag}{\mathrm{diag}}

\theoremstyle{plain}
\newtheorem{thm}{Theorem}[section]

\newtheorem{lem}[thm]{Lemma}
\newtheorem{prop}[thm]{Proposition}
\theoremstyle{definition}
\newtheorem{definition}[thm]{Definition}

\newtheorem{ex}[thm]{Example}
\newtheorem{rem}[thm]{Remark}
\newtheorem{rems}[thm]{Remarks}
\newtheorem{ques}[thm]{Question}

\begin{document} 

\title[Filling families and strong pure infiniteness]
{Filling families and strong pure infiniteness}

\author[Kirchberg and Sierakowski]{E.~Kirchberg 
and A.~Sierakowski}

\date{April 24, 2016} 


\begin{abstract}
We introduce filling families with the matrix 
diagonalization property of 
Definition \ref{def:MatrixDiag}
to extend
and refine  the work by 
R{\o}rdam and the first named author in 
\cite{KirRorOinf}. They allow 
to establish transparent necessary and sufficient 
criteria for strong pure infiniteness of 
crossed products of \cst-algebras by group
actions as obtained in \cite{KirSie1}.
Here we use particular filling families to
improve a result on ``local''
pure infiniteness in \cite{BlanKir2}
and show that the minimal tensor product of 
a strongly purely infinite \cst-algebra and 
a exact \cst-algebra is again strongly 
purely infinite.
Finally we derive with help  of 
suitable filling families an easy
sufficient criterion for the 
strong pure infiniteness of
crossed products 
$A\rtimes_\varphi \N$ by an endomorphism 
$\varphi$ of $A$ 
(\cf  Theorem \ref{thm:endo.cross.spi}).
Our work confirms that the special 
class of nuclear
Cuntz-Pimsner 
algebras constructed in 
\cite{HarKir} consist of strongly purely
infinite \cst-algebras, and thus absorb $\Oinf$ 
tensorially.
\end{abstract}

\maketitle

\tableofcontents

\section{Introduction}
The classification program  of 
G.~Elliott  for nuclear
\cst-algebras \cite{Elliott, Rordam.Classif},
has been an active field of research for
more than 40 years, 
beginning with the classification of AF-algebras
by 
Bratteli \cite{Bratteli}  and Elliott  \cite{Elliott.AFalg}.
This paper focuses
how one might verify when \cst-algebras are 
strongly purely infinite, a property
which is necessary for classification of 
separable nuclear \cst-algebras with 
the help of an ideal system 
equivariant version 
of KK-theory.

In Section \ref{sec2},
following a short Section \ref{sec3} on our 
notation and preliminary results,
we familiarise 
the reader with the notion of
strongly purely infinite \cst-algebras $A$.
One formulation
of this property 
(see Remark \ref{rem:dn.contraction}) 
is that for each  
given pair of positive elements 
$a_1,a_2\in A$, any $c\in A$ and 
$\varepsilon\ge \tau >0$
there exist elements $s_1,s_2\in A$ 
such that
\begin{equation}\label{InEq.general}
\| s_1^* a_1 s_1- a_1 \|<\varepsilon\,, \quad
\| s_2^* a_2 s_2- a_2 \|<\varepsilon\, \quad
\text{and} \quad \| s_1^*cs_2 \| <\tau\,.
\end{equation}
We discus a number of different formulations,
relate the notion of strong pure 
infiniteness to other similar 
notions, and perhaps most importantly
connect it to $\Oinf$ absorption,
classification of 
non-simple \cst-algebras
and previous work in
\cite{DykemaRordam,Kir,
Kir.SPIperm.2003.Preprint,KirAbelProc,
KirRorOinf,TomsWint, Wint1} among others.

In Section \ref{sec:fillFam.Spi}
we introduce the notion 
of a filling family and 
a family with the matrix diagonalization property.
The first notion is roughly speaking 
a intrinsic property encoding a certain
ideal structure for a \cst-algebra
(for a \cst-subalgebra $B\subseteq A$ the map 
$I\mapsto I\cap B$ from
ideals in the \cst-algebra $A$ to ideals in $B$ 
is injective if the positive element in $B$ is 
a filling family for $A$, see Remark 
\ref{rems:filing.familiy}(ii)).
The later notion is a weakening of strong pure 
infiniteness
where we look at solutions of the inequality 
(\ref{InEq.general}) 
but only for a specified family of positive 
elements $a_1,a_2$ in $A$. We prove 
the following result.
\begin{thm}\label{prop:local-spi}
Suppose that $A_+$ contains a filling family
$\cF$ 
(Def.~\ref{def:filling.mdiag.family}),
that has the diagonalization property 
in $A$
(Def.~\ref{def:mdiag.family}).
Then $A$ is strongly purely infinite.
\end{thm}

In Section \ref{sec5}
we develop tools for the verification the
matrix diagonalization property. 
The properties that we study
are of the following type:
Given subsets 
$\cF \subseteq  A_+$, $\Co \subseteq  A$ 
and $\cS \subseteq  A$.
Suppose for each  
given $a_1,a_2\in \cF $, $c\in \Co $ and 
$\varepsilon\ge \tau >0$
there exist elements $s_1,s_2\in \cS $ 
that fulfill
(\ref{InEq.general}).
\begin{ques}
Under which conditions on $\cF $, $\Co $ and 
$\cS $ can the inequalities (\ref{InEq.general}) 
be solved by some $s_1,s_2\in \cS$ for given 
$(a_1,a_2,c,\varepsilon\ge \tau >0)$ with 
$a_1,a_2\in \cF $, but with more general 
elements $c$? 
\end{ques}
We show (as a special case of Lemma 
\ref{lem:cS.control}) that the inequalities 
(\ref{InEq.general}) can be solved for any 
$c$ in the closure of the linear span of 
$\Co$ provided that (i) $\cF$ is invariant 
under 
$\varepsilon$-cut-downs, (ii) $\cS$ is a 
(multiplicative) sub-semigroup of $A$ and 
(iii) for every $s,s_1,s_2\in \cS$, 
$\varphi \in \Cf_c (0,\infty]_+$ and $c\in \Co$
$$\varphi(a_1)s\in \cS,\ \  \varphi(a_2)s\in \cS,
\ \  s_2^*\Co s_1\subseteq  \Co,\ \ 
\textrm{ and } 
\varphi(a_1)c\varphi (a_2)\in \Co.$$
Our results are more general as applications 
require the study of families $\cF $ that are 
not necessarily invariant under 
$\varepsilon$-cut-downs.

In Section \ref{sec6} we consider our first 
application: tensor products. By invoking on 
the work in \cite{BlanKir2} we show the following 
result (where it does not matter which of the 
two algebras is exact):
\begin{thm}
\label{thm:A.ot.B.spi.if.A.spi.B.exact}
The minimal tensor product $A\otimes^{min} B$
of a
\cst-algebra $A$ and an exact \cst-algebra 
$B$ is strongly purely infinite if at least 
one of $A$ or $B$
is a strongly purely infinite \cst-algebra.
\end{thm}

In Section \ref{sec:endo.cross} we consider our 
second application: 
endomorphism crossed products. 
We begin the section by introducing the action
$\sigma\colon \Z \to \mathrm{Aut}(A_e)$ 
associated to $\varphi$,
which is the corresponding
action of the integers 
$\Z$ on the inductive limit
$A_e$ of 
the sequence 
$\xymatrix{A \ar[r]^\varphi & A \ar[r]^\varphi & 
A \ar[r]^\varphi & \cdots }\,$. 
We then require mild sufficient
properties (ND) and (CP) ensuing that 
$A\rtimes_\varphi \N$ 
can be naturally identified with a 
hereditary \cst-subalgebra of 
$A_e\rtimes_{\sigma}\Z$. 
The properly (ND) is automatic when $\varphi$ is 
injective and the second property (CD) ensures 
that the canonical map $A\to A_e$
extends to a strictly continuous 
$^*$-homomorphism
$\cM (A)\to  \cM( A_e)$ of the multiplier 
algebras, \cf{} Lemma \ref{lem:preserve.corners}. 
We prove the following result (as a special case 
of Theorem \ref{thm:endo.cross.spi} where it is 
only required that $b_1,b_2, c$ belongs to a dense 
$\varphi$-invariant \cst-local *-subalgebra of $A$):
\begin{thm}\label{thm1.4}
Let $\varphi$ be an endomorphism of a separable 
\cst-algebra $A$ satisfying properties (ND) and 
(CP). 
Suppose that $\sigma$
is residually properly outer 
(Def.~\ref{def:properly.outer})
and that for every 
$b_1,b_2, c\in A$ and
$\varepsilon>0$ 
there exist 
$k,\,n_1,\,n_2\, \in \N\cup \{0\}$ and 
$s_1,\, s_2\,\in A$ such that
$\| \,  s_j^*\varphi^k(b_j^*b_j)\,s_j-
\varphi^{n_j} (b_j^*b_j)\,\|  \, <\varepsilon$ 
for $j=1,2$ and
$\| \, s_1^*\varphi^k(c)s_2\,\|  \, <\varepsilon$. 
Then $A_e\rtimes_{\sigma}\Z $ and 
$A\rtimes_\varphi \N$ is strongly purely infinite.
\end{thm}

We end by looking at certain class of Cuntz-Pimsner 
algebras. 
More specifically
we look at $^*$-mono\-mor\-phisms 
$h\colon C\hookrightarrow \cM (C)$ of stable 
nuclear separable $\sigma$-unital \cst-algebras 
$C$ to which we associate its canonical 
Hilbert bi-module 
$\mathcal{H}(h,C)$. It turns out that for many 
cases of interest -- 
identifiable in 
terms of conditions on $h$ -- the Cuntz-Pimsner
algebra 
$\mathcal{O}(\mathcal{H}(h,C))$ is strongly 
purely infinite and hence tensorially absorb 
$\Oinf$ (see Remark \ref{rems:HHEK.discussion}).
We prove these result, 
previously 
shown in 
\cite{HarKir}, by identifying these Cuntz-Pimsner
algebras as endomorphism crossed products.

\section{Notation and preliminary results}
\label{sec2}
Let $A_+$ denote the positive elements of $A$,
and
$a_+:=(|a|+a)/2\in A_+$ and
$a_-:=(|a|-a)/2)\in A_+\,$ the positive 
and negative parts of an selfadjoint element
$a$ of $A$, where $|a|:=(a^*a)^{1/2}$. 
If $a\in A_+$, 
then $(a-\varepsilon)_+$, \ie the positive part of
$a-\varepsilon 1\in \cM (A)$, is again in $A_+$ itself.
Here  $\cM (A)$ is the 
\emph{multiplier algebra} of $A$. 
This notation
will be used also for functions $f\colon \R\to \R$,
then \eg 
$(f-\varepsilon)_+(\xi)= 
\max (f(\xi)\, -\varepsilon\,, 0)\,$.
Clearly, 
$\chi ((a-\varepsilon)_+)=(\chi(a) -\varepsilon)_+$ 
for  each
character $\chi$  on the \cst-subalgebra  
$C^*(a)\subseteq  A$  generated by $a$.
This  implies for all $a\in A_+$ and $b\in A$ that
$$\| (a-\varepsilon)_+ \| = ( \| a\| -\varepsilon )_+ 
\quad   \text{and} \quad
\| b-(a-\varepsilon)_+\| \leq \| b-a \| +\varepsilon
\,.$$
A subset $\cF\subseteq  A_+$ is 
invariant under \emph{$\varepsilon$-cut-downs} 
if for each $a\in \cF$ and 
$\varepsilon \in (0,\| a\|)$
we have
$(a-\varepsilon)_+\in \cF$.
The minimal unitalization of $A$
is denoted $\tilde{A}$.
Restriction of a map $f$ to $X$
is denoted $f|X$.
We
let $\Cf_c (0,\infty]_+$ 
denote the set of all non-negative continuous
functions $\varphi$ on $[0,\infty)$ with 
$\varphi | [0,\eta] = 0$ 
for some $\eta \in (0,\infty)$, 
such that $\lim_{t\to \infty} \varphi(t)$ exists.

\begin{rems}\label{rems:b-eps=d*ad}
(i)
Suppose that $\,a,b\in A_+\,$ 
and $\,\varepsilon > 0\,$ 
satisfy  $\| \,  a-b\, \|  \, <\varepsilon\,$. 
Then
$\,(b-\varepsilon)_+\in A\,$
can be decomposed into 
$\,d^*ad=(b-\varepsilon)_+\,$ with some 
\emph{contraction}
$\,d\in A\,$
(\cite[lem.~2.2]{KirRorOinf}).
\smallskip

\noindent 
(ii)
Let $\tau \in [0,\infty)$ and 
$0\leq b\leq a+\tau \cdot 1$ 
(in $\cM  (A)$), 
then for every 
$\varepsilon >\tau$ there is a 
\emph{contraction}
$f\in A$ such that $(b-\varepsilon)_+=f^*a_+f$.
(See
\cite[lem.~2.2]{KirRorOinf} and
\cite[sec.~2.7]{BlanKir2}.)
\smallskip

\noindent
(iii)
Let $a,b,y\in A$, $\delta >0$. 
Let $\varphi$ denote the continuous function 
with compact support in $(0,1]$ given
by $\varphi(t):= \min (1, (2/\delta)(t-\delta/2)_+)$ 
and the $q:=q(y,\delta)\in A$ be the contraction
$q(y,\delta):= \varphi(yy^*)v=v\varphi(y^*y)$, where
$v|y|=y$ is the polar-decomposition 
of $y$ in $A^{**}$. 
Then:
\begin{itemize}
\item[(1)] 
$q^*(yy^*-\mu)_+q=(y^*y-\mu)_+\,$ 
for all $\mu \ge \delta\,$.
\item[(2)]
$0\leq a\leq b$ and $ab=a$  imply   
$a(b-1/2)_+=a/2$.
\end{itemize}
\smallskip

\noindent
(iv)
A  matrix $[b_{k\ell}]\in M_2(A)$ is positive, 
if and only if,
$b_{11},b_{22}\in A_+$, $b_{21}=b_{12}^*$ and 
the transformation
$(b_{11}+1/k)^{-1/2}\, b_{12}\, (b_{22}+1/k)^{-1/2}$ 
is a contraction for every $k\in \N$.
If $[b_{k\ell}]\in M_2(A)_+$, then 
$b_{12}\,=\, \lim_{k\to \infty} \,
b_{11}^{1/2}\,(b_{11}+1/k)^{-1/2}\, 
b_{12}\, (b_{22}+1/k)^{-1/2}\, b_{22}^{1/2}$.

If $[a_{ij}]\in M_n(A)_+$, $n \geq 2$, then the 
$2\times 2$-matrices $[b_{k\ell}]\in M_2(A)$ 
with $b_{11}:=a_{ii}$, $b_{22}:= a_{jj}$ and  
$b_{21}^*=b_{12}:=a_{ij}$
are positive
for each $i\neq j$.
In particular,
$a_{ij}\in \overline{ a_{ii}A a_{jj} }$,
and  
$a_{ij}=\lim _{k\to \infty}  b_i^{(k)} a_{ij} b_j^{(k)}$ 
for the contractions 
$b_j^{(k)}:=(a_{jj}+1/k)^{-1/2} a_{jj}^{1/2}$.

We omit the proofs of (i)--(iv): 
They 
are cited or can be checked easily.
\end{rems}

\section{Strongly purely infinite 
\cst-algebras}\label{sec3}
Recent classification theory
for \cst-algebras in the bootstrap class 
\cf \cite[chp.~9.23]{Blackadar.Ktheory}
extends partly to non-simple 
algebras
 (\cf{}\cite{ref6, ref7, Kir, ref17, ref22, ref24}). 
%%%
%%%
The classification of non-simple nuclear
\cst-algebras 
requires to take in account the structure of the 
primitive ideal spaces. If we classify algebras
with the help of an ideal 
system equivariant version 
of KK-theory, then we can not distinguish
an algebra $A$ from $A \otimes \Oinf$.
This is 
because one can tensor the 
ideal system equivariant 
KK-equivalences with ordinary 
KK-equivalences of nuclear algebras
(\footnote{\, See  
\cite[prop.~2.4(b)]{Skandalis.1988}
for the non-equivariant case,
the equivariant
case is similar.})
and then one can use that  $\Oinf$ is 
KK-equivalent to the complex numbers $\C$
($\Oinf$ and $\C$ have the same K-theory
by \cite[cor.~3.11]{Cuntz:K-Th.in.Annals} 
and are in the bootstrap class, hence are 
KK-equivalent).
Thus, the class of algebras suitable for such 
a classification
contains only nuclear
separable \cst-algebras 
that absorb $\Oinf$ tensorially. 

The requirement $A\ot \Oinf \cong A$ 
looks like a simple criterium, 
but is difficult 
to verify, 
\eg for crossed products. An
intrinsic characterization of 
$\Oinf$ absorbing nuclear
separable 
\cst-algebras 
motivated the following
notion of \emph{strongly purely infinite} algebras:
\begin{definition}\label{def:spi}
A \cst-algebra $A$ is
\emph{strongly purely infinite}
(for short: \textit{s.p.i.}\,)  if,
for every $a, b\in A_+\,$ and  $\varepsilon>0$, 
there exist elements $s,t\in A$ such that
\begin{equation}
\label{InEq.spi}
\| \,  s^*a^2s - a^2 \,\|  \, <\varepsilon\,,\;
\| \,  t^*b^2t - b^2\,\|  \, <\varepsilon
\;\;\textrm{and}\;\; \| \, s^*abt\,\|  \, 
<\varepsilon\, .
\end{equation}
\end{definition}

It was shown in \cite{KirRorOinf} that every 
$\Oinf$ absorbing \cst-algebra is 
strongly purely infinite. 
If the \cst-algebra $A$ is separable, 
nuclear and strongly purely infinite then, 
conversely, $A$ tensorially absorbs $\Oinf$ 
(\cf  \cite{KirRorOinf} for the 
cases of stable or unital algebras, 
and \cite[cor.~8.1]{Kir.SPIperm.2003.Preprint}
for the general case, see also
\cite{TomsWint, Wint1} or 
\cite[prop.~4.4(5), rem.~4.6]{KirAbelProc}
for other proofs of the general case).
There exist strongly purely infinite  non-nuclear 
stable simple separable C*-algebras $A$, 
that are not
isomorphic to $A\ot \Oinf$, 
\cf \cite{DykemaRordam}.
See also \cite[exp.~4.6]{KirRor.2014} for an
example of 
a simple and exact 
crossed product $A\rtimes G$ 
of a type~I  \cst-algebra $A$
by the exact group 
$G:= F_2\times \Z$ safisfying 
that $A\rtimes G$ is strongly 
purely infinite but does not
absorb $\Oinf$.

The nuclearity of the algebra is not a natural
assumption for the study of 
strong pure infiniteness, 
because 
proofs  for 
$\mathrm{KK}$-classification use  
corona algebras
or  asymptotic algebras, 
that are even not exact for
not sub-homogenous algebras  
(\footnote{\, It is because   
$\cL  (\ell_2)$ is a \cst-subquotient
of each not
sub-homogenous  
sub-Stonean algebra.}),
but are still strongly purely infinite
in the sense of 
Definition \ref{def:spi}. 
Fortunately,  multiplier algebras, 
stable coronas and asymptotic algebras of 
strongly purely infinite 
$\sigma$-unital algebras are 
again strongly purely infinite. 

The very basic result for the 
classification program 
is the embedding result for exact algebras into 
strongly purely infinite algebras, \cf \cite{Kir}. 
In this way the notion of 
strongly purely infinite algebras  is 
of  importance for  the classification program. 
This explains our desire to find methods 
and criteria that allow to check
if a given class of (not necessarily simple) 
algebras are 
purely infinite in the \emph{strong} sense of
Definition \ref{def:spi}.

It has been realized 
in an early stage of the
classification 
of -- simple -- pi-sun\  algebras
that many of those algebras are 
stably isomorphic to crossed product
of boundary actions of hyperbolic groups
\cite{AnantharamanDelaroche1997, 
LacaSpi:purelyinf,
JolissaintRobertson}
or as corner-endomorphism cross-product
\cst-algebras 
\cite{Cuntz:On, Cuntz:Symposia38.1982}
and its generalizations.  Therefore it is
likely that criteria for strong pure
infiniteness of crossed products 
can be helpful to detect also
the range of 
KK-classification of non-simple 
\cst-algebras.

\begin{rems}\label{rem:dn.contraction}
(i)
It was shown in the 
proof
of \cite[cor.~7.22]{KirRorOinf} --
but not mentioned in its formulation --
that 
the Definition \ref{def:spi} of 
strong pure infiniteness 
implies that for each 
$a, b\in A_+\,$ and $\,c\in A$
there exist contractions
$s_1,s_2\in A_\omega$
with $s_1a=as_1$, 
$s_2b=bs_2$,
$s_1^* s_1a=a$, 
$s_2^* s_2b=b$
and $s_1^*cs_2=0$.
In particular we have that
our Definition 
\ref{def:spi} 
of strongly purely infinite
\cst-algebras
is equivalent to  
the 
formally stronger requirement, 
that for each 
$a, b\in A_+\,$, $\,c\in A$
and $\varepsilon>0$,
there exist \emph{contractions}  $s,t\in A$
such that
\begin{equation}
\label{InEq.spi.final}
\| \,  s^*as - a\,\|  \, <\varepsilon\,,\;
\| \,  t^*b t - b\,\|  \, <\varepsilon
\;\;\textrm{and}\;\; \| \, s^*c t\,\|  \, 
<\varepsilon\,.
\end{equation}
The proof of 
\cite[cor.~7.22]{KirRorOinf} contains 
some typos    
(\,\footnote{\, Replace in \cite{KirRorOinf}
``$d \in B_\omega$''  
by
``$d\in M_2(A)_\omega$''
and 
``$D$ of $B_\omega$'' 
by 
``$D$ of $M_2(A)_\omega$'' 
in line +14  on page 252,  
``$B_\omega = M_2((B_0)_\omega)$''
by 
``$M_2(A)_\omega = M_2(A_\omega)$'' 
in line -6 on page 252,
and
 ``contractions  in $B_0\subseteq A$''  
by   
``contractions  in $A$''
in line +11 on page 253.}\,).
Compare also the proof of the 
implication (s.p.i.)$\Rightarrow$(I)
in 
\cite[thm.~4.1]{Kir.SPIperm.2003.Preprint}.
\smallskip

\noindent 
(ii)
The proofs  of  \cite[cor.~7.22]{KirRorOinf}
and of 
\cite[thm.~4.1]{Kir.SPIperm.2003.Preprint}
giving contractions $s,t$ in inequalities
\eqref{InEq.spi.final}
both 
 use a fairly  deep local
version of a
``generalized Weyl-von Neumann theorem'' 
\cite[thm.~7.21]{KirRorOinf}. 
But if we apply our Lemma 
\ref{lem:flexible.diag.imply.general}(iii)
to $A$ and  $\cF:=A_+$
then we get 
at least the following  bounds:
\emph{A \cst-algebra $A$ is 
strongly purely infinite in the sense of 
Definition 
\ref{def:spi},  if and only if, for every 
$a_1, a_2 \in A_+$,  $c\in A$
and $\varepsilon \ge \tau>0$,
there  exists $s_1,s_2\in A$ 
that satisfy the inequalities  
\emph{(\ref{InEq.general})}
and have norms 
that satisfy $\| s_j \|^2 
\leq 2\| a_j \|/\varepsilon$.}

\smallskip

\noindent 
(iii)
If we take $a=b$ in the Definition 
\ref{def:spi},
then the
inequalities show that each 
$a^2\in A_+$ is 
properly infinite in $A$, 
see \cite[prop.~5.4]{KirRorOinf}.
Hence  $A$
is purely infinite in sense of 
\cite[def.~4.1]{KirRor1}
by \cite[thm.~4.16]{KirRor1}. In 
general it is an open question whether the 
notions of strong pure 
infiniteness, pure infiniteness and weak pure 
infiniteness coincide or not. When $A$ is 
simple the three properties are equivalent.
We refer 
to \cite{BlanKir2} and \cite{KirRorOinf} for
other special cases where 
weak and strong pure infiniteness coincide.
\smallskip

\noindent 
(iv)
It should be noticed 
that the original definition of J.\ Cuntz 
of purely infinite \cst-algebras in 
\cite{Cuntz:K-Th.in.Annals}
coincides only in some special cases
-- \eg  for simple algebras --  
with the definition
of purely infinite \cst-algebras in 
\cite{KirRor1}. 

Also \cite[thm.~9]{LacaSpi:purelyinf} 
does not show  pure infiniteness 
for crossed products coming from 
\emph{local} boundary
actions \cite{LacaSpi:purelyinf} -- 
even not in the sense of 
\cite[def.~4.1]{KirRor1}.
Both these definitions in 
\cite{Cuntz:K-Th.in.Annals, KirRor1} are 
still not suitable for the classification 
in general
 --
except in combination
with other 
assumptions, like \eg  tensorial  absorption of the 
Jiang-Su algebra $\mathcal{Z}$.
\end{rems}

\section{Filling families and 
strong pure infiniteness}%1362
\label{sec:fillFam.Spi}
A suitable algebraic theory for 
invariants of strongly purely
infinite \cst-algebras is not in sight, 
different to the property
of pure infiniteness of \cst-algebras 
$A$ that is
equivalent  to $2[a]=[a]$ in the Cuntz 
semigroup
for every $a\in (A\ot \cK)_+$. 
The absence of a reasonable algebraic
description 
forces us to develop new methods to detect 
and describe strong pure infiniteness.

Here we introduce two new concepts: 
Firstly we work with 
the idea of a \emph{filling} family 
$\cF\subseteq A_+$, \cf Definition 
\ref{def:filling.mdiag.family}.
Secondly we introduce the notion of 
a family $\cF\subseteq A_+$ with the 
\emph{matrix diagonalization property} as 
a refinement of the notion of 
the matrix diagonalzation property introduced
in \cite[def.~5.5]{KirRorOinf} 
(\cf Definition \ref{def:MatrixDiag}).
We say ``family''  because we use the 
elements of $\cF$ mainly to form a 
family of selfadjoint $n\times n$-matrices 
with diagonal entries
from $\cF$ for $n=2,3,\ldots$ -- 
together with certain restrictions
on the off-diagonal entries.

Before defining a filling family 
we need a lemma.
Notice that one can replace 
in part (ii) of the following 
Lemma \ref{lem:filling.family}
\emph{primitive} ideals by 
\emph{all closed} ideals $I$ with 
$D\not\subseteq I$,
because every closed ideal is the 
intersection of primitive ideals.
\begin{lem}\label{lem:filling.family}
Let $\cF $  be a subset of  $A_+$.
The following 
properties of $\cF $ are equivalent:
\begin{itemize}
\item[(i)]
For every $a,b,c\in A$ with 
$0\leq a\leq b\leq c\leq 1$,
with 
$ab=a\not=0$ and  $bc=b$, there exists
$z_1,z_2,\ldots, z_n\in A$ and $d\in A$ with
$z_j(z_j)^*\in \cF $, such that
$ec=e$
and $d^*ed=a$ for 
$e:= z_1^* z_1+\ldots+ z_n^* z_n$.
\item[(ii)] 
For every hereditary \cst--subalgebra
$D$ of $A$ and every primitive ideal
$I$ of $A$ with $D\not\subseteq I$ there
exist $f\in \cF $ and $z\in A$
with $z^*z\in D$ and  $z z^* =f  \not\in I$.
\end{itemize}
\end{lem}

\pf (i)$\Rightarrow$(ii):
Since $D\not\subseteq I$, 
there exists
$g\in D_+$ with $\| g\|=\| g+I \|=3$.
Let $a:=(g-2)_+$,
$b:=\, (g-1)_+ - (g-2)_+\,$ and 
$c:=\, g-(g-1)_+\,$. 
Then $a\not\in I$,
$0\leq a\leq b\leq c\leq 1$, $ab=a$, $bc=c$,
$g=c+b+a$ and $\|c+b\| \leq 2$. 
By (i),
we find $z_1,z_2,\ldots,z_n\in A$ and $d\in A$ 
with
$z_j(z_j)^*\in \cF $, such that
$e:=z_1^*z_1+\ldots+z_n^*z_n$ satisfies 
$ec=e$ and $d^*ed=a$.
It follows that 
$e\in cAc\subseteq D$ and $e\not\in I$.
Hence $(z_j)^*z_j\in D\backslash I$ for some 
$j\in \{ 1 ,\ldots, n\}$.
Then $z:=z_j$ and 
$f:=z_j(z_j)^*\in \cF $ satisfy (ii). 

(ii)$\Rightarrow$(i): 
Suppose that $a,b,c\in A$ with 
$0\leq a\leq b\leq c\leq 1$,
$ab=a\not=0$, $bc=b$ are given. 
Let  $D:=\overline{bAb}$ and let
$\mathcal{Z}$ denote the set of
all $z\in A$ with $z z^*\in \cF $ and $z^*z\in D$.
It follows that $z^*z c = z^*z$  and
$zu\in \mathcal{Z}$ for $z\in \mathcal{Z}$
and every unitary 
$u$ in the minimal unitization 
$\widetilde{D}$ of $D$. 
Consider  the set $M$  of
$d \in D_+$ with the property that there are
$z_1,\ldots, z_n\in \mathcal{Z}$ and a
$\rho \in (0,\infty)$ with
$$0\leq d\leq 
\rho \cdot (z_1^*z_1+ \ldots +z_n^*z_n)\,.$$
Clearly,  $M$ is a (not necessarily closed) 
hereditary convex cone in $D_+$,
and  $u^*M u\subseteq M $ for every unitary 
$u$ in $\widetilde{D}$. 
Arguments in the  proof of 
\cite[thm.~1.5.2]{Ped.book} 
show that the set 
$L(M):= \{ d\in D\,;\,\, d^*d \in M\,\}$ 
is a (not necessarily closed)
left ideal of $D$. Since each element of an
unital \cst-algebra is the linear
combination of  unitaries  in this algebra, 
we get that $d x \in L(M)$ for 
each $x\in \widetilde{D}$, $d\in L(M)$. 
In particular, $L(M)$ is a two-sided ideal
of $D$.
It follows, that the closure 
$J:= \overline{ L( M ) }$  of the two-sided ideal 
$L( M )$ is an  ideal of $D$.
The closure $K:= \overline{\,M\,}$ of 
$M$ is again a hereditary convex cone in 
$D_+$, with $L(K)=J$. 
Thus, $K$ is the positive part 
$J_+$ of the closed left ideal $J$ of $D$.

Since $D$ is hereditary, 
the closed linear span $I$ of $AJA$ 
is a closed ideal $I$ of $A$ with $J=D\cap I$, 
see 
\cite[II.5.3.5]{Blackadar.CstAlg}.

The property (ii) implies  $D\subseteq I$,  
because $D\cap I$ 
contains \emph{all} elements $z^*z$ with 
$z\in \mathcal{Z}$ 
(since $z^*z\in M$ for $z\in \mathcal{Z}$). 
We obtain that $D_+=J_+=\overline{\,M\,}$.
In particular, $b\in \overline{\,M\,}$, and
we find an element $g\in M$ with 
$b\leq g + 1/4$.
By definition of $M$ there are 
$z_1,\ldots,z_n\in \mathcal{Z}$
and $\rho\in (0,\infty)$ with
$g\leq \rho e$ for
 $e:=z_1^*z_1+ \ldots +z_n^*z_n$.
 Then $e\in D_+$ and $ec=e$
 by definitions of $D$ and of $\mathcal{Z}$.
There is
a contraction $d_1\in D$ with  
$(b-1/2)_+=\rho d_1^*ed_1$ by 
Remark \ref{rems:b-eps=d*ad}(ii).
It follows $a=d^*ed$ for 
$d:=\sqrt{2\rho}\cdot d_1 a^{1/2}$,
because  $a(b-1/2)_+=a/2$.
by Remark \ref{rems:b-eps=d*ad}(iii,2).
\qed

\begin{definition}\label{def:filling.mdiag.family}
Let $\cF $  be a subset of  $A_+$.
The set $\cF $ is a 
\emph{filling family} for $A$, if
$\cF $ satisfies the equivalent conditions (i) 
and (ii)
of Lemma \ref{lem:filling.family}.
\end{definition}

\begin{rems}
(i)
Lemma \ref{lem:filling.family}(ii) shows
that $\cF:=A_+$ is a filling 
family for $A$.
	
\noindent 
(ii)
We warn the reader that our fundamental
notion of ``filling families'' does not require
an existence of some type of actions,
e.g.\ group actions.
In particular, they have nothing in 
common with the notion of $n$-filling
actions in \cite{JolissaintRobertson}.
\smallskip

\noindent 
(iii)
\label{rem:simpler.if.cF.has.cut.downs}\,
If $\cF\subseteq  A_+$ is 
\emph{invariant under 
$\varepsilon$-cut-downs},
\ie if $(a-\varepsilon)_+\in \cF$
for each $a\in \cF$ and 
$\varepsilon \in (0,\| a\|)$,
then we can replace the Murray--von-Neumann
equivalence $z^*z \approx_{MvN} f $ in 
Lemma \ref{lem:filling.family}(ii)
by  Cuntz equivalence (denoted by $\sim$)
in Definition \ref{def:filling.mdiag.family}(i):\\
\emph{For every hereditary \cst--subalgebra
$D$ of $A$ and every primitive ideal
$I$ of $A$ with $D\not\subseteq I$ there
exist $f\in \cF \backslash I$ and $g\in D$
with $g\sim f$}.
\end{rems}

\begin{rems}\label{rems:filing.familiy}
It is in general not easy 
to check if a family 
$\cF \subseteq  B_+$ is 
filling for $B$ or not. We list some cases 
$\cF \subseteq  A_+\subseteq  B$,
where $A\subseteq B$ are \cst-algebras:
\begin{itemize}
\item[(i)]  If $A=\Cf_0(X)$, then 
$\cF\subseteq  A_+$ is filling,
if and only if, the supports of the functions 
$f\in \cF$ 
build a base
of the topology of $X$. 
\item[(ii)] If 
$\cF:=A_+\subseteq B$ is filling for $B$, then
the map 
$I\in \mathcal{I}(B)\mapsto
 I\cap A \in  \mathcal{I}(A)$ 
is injective, \ie $A$ separates the 
closed ideals
of $B$, \cf \cite{Siera2010}.
\item[(iii)] Let $R$ denote a nuclear, 
separable, 
simple and finite 
\cst-algebra
such that $B:= M_2(R)$ is properly infinite 
(\cf R{\o}rdam 
\cite{Rordam.Example}).
Thus, there exist a unital *-homomorphism 
$\varphi\colon \Oinf \to B$.
Consider the image $A:= \varphi(\Oinf)$
and let $\cF:=A_+$.
Then $\cF$ separates the ideals 
$\{ 0\}$ and  $B$  of $B$, 
but $\cF$  is not filling for $B$.
\item[(iv)] If $D\not= \C$ is a simple, unital, and 
\emph{stably finite} \cst-algebra, $X$ is a 
locally compact Hausdorff space, 
$B:=\Cf_0(X)\otimes D$ and 
$A:=\Cf_0(X) \otimes 1 \subseteq B$,
then $\cF := A_+$  separates the ideals of $B$, 
but  is not filling for $B$.
\item[(v)] Let 
$A:=\mathbb{C}\oplus \mathbb{C}$, 
$\sigma(u,v):=(v,u)$, 
$\cF:=\{ (1,0) \}\subseteq A_+$, and 
$B:=A\rtimes _\sigma \Z_2\cong M_2(\C )$.
Then $\cF$ and $A_+$ are filling for $B$, but 
$\cF$ is not filling for $A$.
\item[(vi)] The example in \cite[p.123]{ArchSpiel}
has the property that 
$\cK \ot A_{1-\theta} \cong \cK \rtimes \Z^2$
is simple and stably finite and 
$\cF:=\cK_+ \subseteq \cK \rtimes \Z^2$
is not filling for  $\cK \rtimes \Z^2$.
\end{itemize}

All above examples are refereed or easily verified.
For the convenience of the reader we give some
hint for part (vi):  The only 
primitive ideal  of $\cK \rtimes \Z^2$ is 
$\{0\}$. Take  a minimal projection 
$p\in \cK_+\subset \cK \rtimes \Z^2$. Then we
find a non-zero projection $q\in p(\cK \rtimes \Z^2)p$
that is not MvN-equivalent to 
$p$. Let $D:=q(\cK \rtimes \Z^2)q$.
Then there is no non-zero contraction 
$z\in \cK \rtimes \Z^2$ 
with $z^*z\in \cK$ and $zz^*\in D$, because $p$
is not  infinite.
\end{rems}

\begin{lem}
\label{lem:F.fill.A.A+.fill.B.Then.F.fill.B}
Suppose that $A\subseteq B$ is a 
\cst-subalgebra of $B$
and
$\cF\subseteq  A_+$ is a subset of $A_+$.
If $\cF$ is filling for 
$A$, and $A_+$ is filling for $B$, then $\cF$
is a filling family for $B$.
\end{lem}
\pf
Let $D\subseteq  B$ hereditary, 
$I\subseteq  B$ closed ideal with 
$D\not\subseteq I$.
By assumption, there is  
$z\in B$ with $z^*z\in D$, 
$zz^*\not\in I$ and $zz^*\in A_+$.
Let $E$ denote the hereditary \cst-subalgebra 
of $A$ generated by
$zz^*$, \ie $E:=\overline{\, zz^*Azz^*\,}$. 
Since $zz^*\not\in I$, 
the algebra $E$ is not contained in the ideal
 $J:=A\cap I$  of $A$. 
 By assumption, there exists 
 $y\in A$ with $y y^*\in \cF$,  
$y^*y\in E\subseteq  A$ and $y y^* \not\in J$.
Let $v(z^*z)^{1/2} = z$ 
denote the polar decomposition of 
$z$ in $B^{**}$.
Then $x:=y v\in B^{**}$ satisfies 
$x\in B$, because $y^*y\in \overline{z z^*Bz z^*}$.
Moreover, $xx^*\not\in I$,
$xx^*\in \cF$ and $x^*x\in v^*Ev\subseteq D$:
To see $xx^*\in \cF$, notice that
$z z^*vv^*=z z^*$, hence for all $e\in E$,
$e vv^*=e$.
Since $y^*y\in E$, we get $\| y vv^*-y\|^2=0$,
so $xx^*=y vv^*y^*=y y^*\in \cF$.
\qed

\begin{definition}\label{def:MatrixDiag}
Let $\cS\subseteq A$ 
be a multiplicative sub-semigroup 
of a \cst-algebra $A$ and $\Co \subseteq A$
a subset of $A$.
An $n$-tuple $(a_1,\ldots,a_n)$ of
positive elements in $A$ has the 
\emph{matrix diagonalization property 
with respect to $\cS$ and $\Co $},
if for every 
$[a_{ij}]\in M_n(A)_+$ with $a_{jj}=a_j$ and 
$a_{ij}\in \Co $
(for $i\not=j$)
and $\varepsilon_j >0, \tau>0$ 
there are elements 
$s_1\,,\ldots , s_n\, \in \cS$ 
with
\begin{equation}
\label{InEq.DiagProp}
\| s_j ^* a_{jj}\, s_j - a_{jj} \| 
< \varepsilon_j\,, \quad
\text{and} \quad 
\| s_i ^*a_{i j}s_j \| < \tau \,\,\,  
\text{for} \,\, i\not= j\,.
\end{equation}
If $\cS=\Co =A$ then this is the 
\emph{matrix diagonalization property} 
of $(a_1,\ldots,a_n)$ as defined in 
\cite[def.~5.5]{KirRorOinf},
and we say that $(a_1,\ldots,a_n)$ 
has matrix diagonalization (in $A$).
\end{definition}
\begin{definition}\label{def:mdiag.family}
Let $\cF $  be a subset of  $A_+$.
The family $\cF $ has the 
\emph{(matrix) diagonalization property}
(in $A$)  if each finite sequence 
$a_1,\ldots, a_n\in \cF $
has the matrix diagonalization property 
(in $A$) of 
Definition \ref{def:MatrixDiag}.
\end{definition}

\begin{rems}
(i)
By Remark \ref{rems:b-eps=d*ad}(i) 
and a preceding inequality, it follows that 
the $n$-tuple $(a_1,\ldots,a_n)$ 
has the matrix diagonalization 
with respect to $A$ and $\Co $
if, and only if, for each 
$[a_{ij}]\in M_n(A)_+$ with 
$a_{jj}=a_j$ and $a_{ij}\in \Co $
(for $i\not=j$)
and $\varepsilon_j >0, \tau>0$ there are
elements $s_1,\ldots,s_n\in A$ 
that satisfy
the equations and inequalities
\begin{equation}
\label{InEq.DiagProp.cS=A}
s_j ^* a_{j j}\, s_j=(a_{jj}-\varepsilon_j)_+\,,  
\quad \text{and} \quad 
\| s_i ^*a_{i j}s_j \| < \tau \,\,\,  
\text{for} \,\, i\not= j\,.
\end{equation}
\smallskip

\noindent 
(ii)
If we replace
the $\varepsilon_j$ and $\tau$ in  
inequalities (\ref{InEq.DiagProp})  by 
$\varepsilon:=
\min (\,\varepsilon_1\,,\, \ldots 
\,,\, \varepsilon_n\, , \,\tau\,)$,
then this new definition is the same as 
Definition  \ref{def:MatrixDiag} with
$\varepsilon_1=\cdots=\varepsilon_n
=\tau=\varepsilon$.
But the latter
is an equivalent formulation of Definition  
\ref{def:MatrixDiag}. 
\smallskip

\noindent 
(iii)
The following is again 
equivalent to 
the
matrix diagonalization property:
The $n$-tuple  
$(a_1,\ldots, a_n)$  has the  matrix
diagonalization property with respect to $\cS$
and $\Co $, if and only if,
for each positive matrix $[a_{ij}]\in M_n(A)$
with diagonal entries $a_{jj}=a_j$ and 
$a_{ij}\in \Co $
(for $i\not=j$), there exists a sequences 
$s^{(k)}\in M_n(A)$,  $k=1,2,\ldots$,
of diagonal matrices 
$s^{(k)}=\diag ( s_1^{(k)},\cdots, s_n^{(k)})$
with $s_j^{(k)}\in \cS$, such that
$$
\lim_{k\to \infty}  \,\,
\|\, (s^{(k)})^*\, [\,a_{i j}\, ] \, s^{(k)}  - 
\diag (a_1,\ldots , a_n) \,\| \,
= 0
\,.$$
\smallskip

\noindent 
(iv)
It is important for our 
applications to find an estimate of 
$\max_j \| s_j\|^2$
depending only on 
$
\min (\,\varepsilon_1\,,\, \ldots \,,\, 
\varepsilon_n\, )
$
that does not depend 
on 
$\{a_{ij} \,;\,\, j\not=i\}$
or on $\tau >0$.
Therefore, we often use (starting 
from proof of Lemma 
\ref{lem:diag.property.MvN.equivalent}) the 
equivalent formulation
of  Definition \ref{def:MatrixDiag}
with values
$\varepsilon_j:=\varepsilon>0$ and 
independent $\tau >0$, considering 
inequalities
\begin{equation}
\label{InEq.DiagProp2}
\| s_j ^* a_{jj}\, s_j - a_{jj} \| 
< \varepsilon\,, \quad
\text{and} \quad 
\| s_i ^*a_{i j}s_j \| < \tau \,\,\,  
\text{for} \,\, i\not= j\,.
\end{equation}
\end{rems}

\begin{lem}
\label{lem:diag.property.MvN.equivalent}
Let $z_1,\ldots, z_n\in A$ such that
$(z_1^*z_1,\ldots, z_n^* z_n)$ 
has the matrix diagonalization property in $A$.
\begin{itemize}
\item[(i)]
If $1\leq k< n\,$,
$\,e:=z_1^*z_1\,+ \cdots +\, z_k^*z_k\,$
and 
$\,f:=z_{k+1}^*z_{k+1}\,+ \cdots + \,z_n^*z_n$,
then $(e,f)$ 
has the matrix diagonalization property.
\item[(ii)]
The n-tuple $(z_1z_1^*,\ldots, z_n z_n^*)$
has the matrix diagonalization property.
\end{itemize}
\end{lem}

\pf (i): Follows from 
\cite[lem.~5.9]{KirRorOinf}.

(ii):\,
Let $[a_{ij}] \in M_n(A)_+$ with 
$a_{jj}=z_j z_j^*$.
By Remark \ref{rems:b-eps=d*ad}(iv),  
$a_{ij}\in \overline{a_{ii} A a_{jj}}$
and $a_{ij}=
\lim_{k\to \infty} b_i^{(k)}a_{ij}b_j^{(k)}$
for the contractions 
$b_j^{(k)}:
=(a_{jj}+1/k)^{-1/2} a_{jj}^{1/2}\ge 0$.
Consider   the polar decompositions 
$z_j^*= v_j |z_j^*|=v_j a_{jj}^{1/2}$
of the $z_j^*$ in  $A^{**}$.
Then 
$v_j b_j^{(k)}= z_j^*(a_{jj}+1/k)^{-1/2}\in A$,
$v_j z_j z_j^*= z_j^*z_j v_j$, 
$z_j^*z_j=v_j a_{jj}v_j^*$
and  $z_j z_j^*=v_j^* (z_j^*z_j) v_j$.
It follows $v_i a_{ij} v_j^*\in A$ 
and $v_i^*v_i a_{ij} v_j^*v_j= a_{ij}$  for  
$i,j= 1,\ldots ,n$.
The diagonal matrix  
$V:= \diag (v_1, \ldots ,v_n)$
is a partial isometry in $M_n(A^{**})$, 
the matrix 
$[ c_{ij} ]  : = V [a_{ij}]V^*$ is a 
positive matrix
in $M_n(A)$ with diagonal entries 
$c_{jj}=z_j^*z_j$,
and $V^*[c_{ij} ] V= [a_{ij}]$.

Let $\varepsilon\ge \tau >0$. 
By assumption,
there are $e_j\in A$ with 
$$\| e_j^* z_j^*z_j e_j - z_j ^*z_j \|= 
\| e_j^*v_j a_{jj} v_j^* e_j - v_j a_{jj} v_j^*\| 
< \varepsilon/2
%$
\quad \text{and} \quad 
%$
\| e_i^* v_i a_{ij}v_j ^* e_j  \| < \tau/2$$
for $ i,j =1,\ldots,n$ and $i\not=j$.
Let $\delta:= \tau / (2+2(\max_j \| e_j\|)^2)$.
We find $k\in \N$ such that
$\| a_{ij}- f_ia_{ij}f_j \| <\delta$ for
$f_j:=b_j^{(k)}$.
Let $s_j:=(v_j f_j)^*e_j(v_j f_j)$.
Since  $v_j f_j=v_j b_j^{(k)}\in A$,
we get that $s_j\in A$. 
This $s_j$ satisfy
$\| \,s_j^*a_{j j}s_j - 
a_{jj}\, \| \,\,
< \,\, 
\delta\|   e_j\|^2  +  \varepsilon/2 +\delta 
\leq \varepsilon$
and 
$\| s_i^*a_{ij}s_j \| \, < \, 
\delta \| e_j\|^2 + \tau/2\leq \tau$ 
giving (\ref{InEq.DiagProp2}).
\qed

\begin{lem}\label{lem:lem.approx.diag}
Let $a,b\in A_+\,$. Suppose that, for each 
$\varepsilon\in (\, 0\,,\, 
\min(\,\| a \|, \| b\|)/4\,)$,
there exist $e,f\in A_+\,$ and  $d_1,d_2\in A$
such that
\begin{itemize}
\item[(i)] $d_1^*ed_1=(a-3\varepsilon)_+$ and 
$d_2^*fd_2=(b-3\varepsilon)_+\,$,
\item[(ii)] 
$\varepsilon e =\,(a-(a-\varepsilon)_+)e\,$,  
$\,\varepsilon f = (b-(b-\varepsilon)_+)f\,$, and 
\item[(iii)] $(e,f)$ has the 
matrix diagonalization property.
\end{itemize}
Then $(a,b)$ 
has the matrix diagonalization property.
\end{lem}

\pf  Let $[a_{ij}]\in M_2(A)_+$ with $a_{11}:=a$ 
and 
$a_{22}:=b$, $\varepsilon>0$
and $\tau>0$.
We show that there exists  $v_1,v_2\in A$ 
such that
$ v_1^*av_1= (a-4\varepsilon)_+$, 
$v_2^*bv_2=(b-4\varepsilon)_+$,
and $\| v_1^*a_{12} v_2 \| < \tau$. 
If
$4\varepsilon\ge \min(\| a \|, \| b\|)$ let
$v_1:=\lambda(a)$  and  $v_2:= \lambda(b)$
with 
$\lambda(t):= 
t^{-1/2}\cdot (t-4\varepsilon)_+^{1/2}$.
If 
$\varepsilon<\min(\| a \|, \| b\|)/4\,$,
let $e,f, d_1,d_2\in A$ be 
elements with the properties in (i)--(iii).
We define continuous functions  
$\psi$ and $\varphi$ 
on $[0,\infty)$
by
$\psi(t):= \min(1,\varepsilon^{-1} t)$, 
$\varphi(t):= \varepsilon ^{-1} t$ for 
$t\leq \varepsilon $ 
and
$\varphi(t):=t^{-1} \varepsilon $ for  
$t>\varepsilon$.
Notice that  $\varphi(t)t= \varepsilon \psi(t)^2$.

Put 
$c_1:=\varepsilon ^{-1}(a-(a-\varepsilon)_+)=
\psi (a)$,  
$c_2:=\varepsilon ^{-1}(b-(b-\varepsilon)_+)=
\psi (b)$.
Then $e=c_1e$, $f=c_2f$,  
$a\varphi (a)=\varepsilon c_1^2$
and $b\varphi (b)=\varepsilon c_2^2$.
Since the elements are all positive, we get  that 
$e$ commutes with
$c_1$ and that $f$ commutes with $c_2$. 
It follows that
$e=e^{1/2}c_1^2 e^{1/2}$ and 
$f=f^{1/2}c_2^2 f^{1/2}$.
We let  
$g_1:= 
\varepsilon^{- 1/2} \varphi(a)^{1/2} e^{1/2}$
and 
$g_2:=
\varepsilon^{- 1/2} \varphi(b)^{1/2} f^{1/2}$.
The $2\times 2$-matrix 
$[b_{ij}]: 
= \diag(g_1,g_2)^*[a_{ij}]\diag(g_1,g_2)$
is positive and has entries 
$b_{11}=g_1^*ag_1=e$,
$b_{22}=g_2^*bg_2=f$, and  
$b_{21}^*=b_{12}=g_1^*a_{12}\,g_2$.

Let 
$\gamma:=
\max \bigl(\,1\,,\, \| d_1\|^2\,,
\, \| d_2\|^2\,\bigr)$ and 
$\delta:= \min(\varepsilon, \tau) / \gamma$.
The  diagonalization property of $(e,f)$
gives $S_1,S_2\in A$ with 
$\| S_1^*eS_1-e \| < \delta$,  
$\| S_2^*fS_2-f \| < \delta$, 
and $\| S_1^* b_{12}\, S_2 \| <\delta$.

It follows that $T_j:=S_j d_j$ ($j=1,2$) satisfy
$$\| T_1^*eT_1 - (a-3\varepsilon)_+ \| 
< \delta \gamma\leq \varepsilon\,, \,\,\,
\| T_2^*fT_2 - (b-3\varepsilon)_+ \| <  
\varepsilon \,\,\,
\text{and} \,\,\, \|  T_1^* b_{12} T_2\| 
<\delta \gamma \leq \tau
\,.$$

Thus, $h_j:= g_jT_j$ satisfies
$\| h_1^*ah_1 - (a-3\varepsilon)_+ \| 
< \varepsilon$,
$\| h_2^*bh_2 - (b-3\varepsilon)_+ \| 
< \varepsilon$
and $\|  h_1^* a_{12} h_2\| < \tau$.
Use Remark \ref{rems:b-eps=d*ad}(i) to 
get the desired
$v_i:=h_i q_i$ with suitable 
contractions $q_i\in A$.\qed

\textbf{Proof of 
Theorem \ref{prop:local-spi}}:\,
Let $a,b\in A_+\backslash \{ 0\}$.
We show that $(a,b)$ 
has the matrix diagonalization
property. 
This applies in particular to the positive matrix
$[a^{1/2},b^{1/2}] ^\top [a^{1/2},b^{1/2}]
\in M_2(A)_+$
and proves that $A$ is 
strongly purely infinite
in the sense of Definition \ref{def:spi}.

Let 
$\varepsilon \in (0, \min(\|a\|,\| b\|)/4)$,
and $\gamma:=(\| a \| - 3\varepsilon)^{-1}\,$.
We show  the existence of $e,f; d_1,d_2\in A$
that satisfy the conditions (i)--(iii)  of 
Lemma \ref{lem:lem.approx.diag}.

By Lemma \ref{lem:filling.family}(i),  we find 
$y_1,\ldots, y_m, d_1\in A$ satisfying
$e \left(a -(a-\varepsilon)_+\right)=
\varepsilon  e$,
$y_i(y_i)^*\in \cF $ and $
d_1^*e d_1=
(a-3\varepsilon)_+\,
$ for $e:= \sum_{i=1}^m  y_i^*y_i$,
because we can apply  
Lemma \ref{lem:filling.family}(i)
to the elements 
$\gamma (a-3\varepsilon)_+\,$,
$\,\varepsilon^{-1} ((a-\varepsilon)_+ -
(a-2\varepsilon)_+)\,$
and 
$\,\varepsilon^{-1} 
(a - (a-\varepsilon)_+)
$
in place of the elements $a\leq b \leq c$
in \ref{lem:filling.family}(i).

In the same way one can see, that
Lemma
\ref{lem:filling.family}(i) gives 
elements  $z_1,\ldots,z_n, d_2\in A$ 
such that
$f\left( b -(b-\varepsilon )_+\right)=
\varepsilon  f$,
$z_j(z_j)^*\in \cF $ and 
$d_2^*f d_2=
(b-3\varepsilon)_+$ for $f:=\sum_j z_j^*z_j$.

Since the sequence 
$(y_1y_1^*,\ldots, y_my_m^*, z_1 z_1^*,
\ldots, z_n z_n^*)$
has the matrix diagonalization property 
(by  assumptions on $\cF$)
the Lemma  
\ref{lem:diag.property.MvN.equivalent} 
applies and shows that the sequences
$(y_1^*y_1,\ldots, y_m^*y_m, z_1^*z_1,
\ldots, z_n^*z_n)$
and $(e,f)$ both
have the matrix diagonalization property.
Thus the elements $e,f\in A_+$ and 
$d_1,d_2\in A$
satisfy the conditions (i)--(iii) 
of Lemma \ref{lem:lem.approx.diag}, and
$(a,b)$ has the matrix diagonalization property by  
Lemma \ref{lem:lem.approx.diag}.
\qed
\vfill

\section{Verification of the matrix 
diagonalization}\label{sec5}
Given subsets 
$\cF \subseteq  A_+$, $\Co \subseteq  A$ 
and $\cS \subseteq  A$. In this section we 
study questions related to the 
verification the matrix
diagonalization property with respect to $\cS$
and $\Co $ for (finite) tuples of elements in $\cF$.
We study questions 
of the following type:

(Q1) \emph{Under which conditions on 
 $\cF $, does it follow that
$\cF$ has the matrix diagonalization 
property?}

(Q2) \emph{Under which  conditions on 
 $\cF $, $\Co $ and $\cS $ can
the inequalities} (\ref{InEq.general}) 
\emph{be solved by some $s_1,s_2\in \cS$
for each given  
$(a_1,a_2,c,\varepsilon\ge \tau >0)$
with $a_1,a_2\in \cF $, and
$c\in\overline{\mathrm{span}(\Co)}$}?

An example of a possible condition for a 
positive answer
to (Q1) is
that $\cF $ is invariant under 
\emph{$\varepsilon$-cut-downs}, 
\ie{}that
for each $a\in \cF$ and 
$\varepsilon \in (0,\| a\|)$
we have
$(a-\varepsilon)_+\in \cF$ (\cf Lemma 
\ref{lem:flexible.diag.imply.general}).
The answer to the second question has to do with
interplay of  $\cF $, $\Co $ and $\cS$.
This means that we have to require additional
suitable conditions, \eg that
$\cS ^*\cdot \Co\cdot \cS \subseteq  \Co$.

We need this generalization because
our applications are concerned with
families  $\cF $ that are not invariant under 
\emph{$\varepsilon$-cut-downs}, \ie  operations
$a\mapsto (a-\varepsilon)_+$ for $a\in \cF $ and 
$\varepsilon \in (0, \|a\|)$. An example is the
proof of
Theorem \ref{thm:A.ot.B.spi.if.A.spi.B.exact}.
It uses the following
Lemma \ref{lem:from.2.diag.to.matrix.diag} 
that we could not directly deduce from 
\cite{KirRorOinf}. We start by a 
definition allowing us to better control the 
matrix diagonalization property:

\begin{definition}\label{def5.1}
A $n$-tuple $(a_1,\ldots,a_n)$
of positive elements in $A$
has \emph{controlled} matrix diagonalization property
with respect to $\cS$ and $\Co $
if there is a non-decreasing  controlling function 
$$
(0,\infty)\ni t\mapsto D_n(t) := 
D_n(t; a_1,\ldots, a_n)
\in [1,\infty)
$$
such that for every 
$[a_{ij}]\in M_n(A)_+$ with $a_{jj}=a_j$ and 
$a_{ij}\in \Co $
(for $i\not=j$)
and $\varepsilon_j >0, \tau>0$ there are
$s_1,\ldots,s_n\in \cS$ that 
satisfy the inequalities  
(\ref{InEq.DiagProp})
and have norms 
that satisfy 
$$
\|s_j \|^2 \leq  
D_n(1/\min(\varepsilon_1,\ldots, 
\varepsilon_n))\,.$$
If $\cS=\Co=A$ we say $(a_1,\dots,a_n)$ has 
controlled matrix diagonalization (in $A$).
\end{definition}
 
The following lemma in parts reduces the 
problem of considering arbitrary $n$-tuples to 
$2$-tuples. We say ``in parts'' because the 
assumptions in Lemma 
\ref{lem:from.2.diag.to.matrix.diag} involve 
matrices in $M_2(A)$ that are not necessarily 
positive. This difficulty is solved in Lemma 
\ref{lem:from.2.diag.Prop.to.Estimate}.

\begin{lem}\label{lem:from.2.diag.to.matrix.diag}
Let $\cF $  be a subset of  $A_+$. 
Suppose that for any  
given $a_1,a_2\in \cF $, 
there exists a non-decreasing 
function $t\mapsto D(t;\, a_1,a_2) <\infty$ 
such that for each $c\in a_1^{1/2}A\,a_2^{1/2}$ 
and $\varepsilon\ge \tau >0$
there exist $s_1,s_2\in A $ that fulfill
\emph{(\ref{InEq.general})}
and 
$\|s_j\|^2 \leq D(1/\varepsilon\,;\,\, a_1,a_2)$.
Then any $n$-tuple of elements in $\cF$ has 
the controlled matrix diagonalization 
in $A$.
\end{lem}

\pf
We can suppose that all the functions 
$t\mapsto D(t;\, a_1,a_2)< \infty$ 
satisfy $D(t;\, a_1,a_2)\ge 1$ 
for all $t\in (0,\infty)$,
upon replacing $D(t;\, a_1,a_2)$  by 
$\widetilde{D}(t;\, a_1,a_2):= 
\max (1, D(t;\, a_1,a_2))$.

Let $a_1,\ldots, a_{n+1}\in \cF $ and 
let $[a_{j k}]$ be a positive matrix 
in $M_{n+1}(A)$ with diagonal 
entries $a_{j j}=a_j$.

We proceed by induction over $n\ge 1$, and
prove each $n$-tuple of elements in $\cF$ has 
the controlled matrix diagonalization. 
It suffice to prove the existence of a 
controlling function
$t\mapsto D_{n+1}(t)=
D_{n+1}(t; a_1,\ldots,a_n) < \infty$ 
with the property
that, for every  $\varepsilon \ge \tau >0$,
there exists  $s_1, \ldots, 
s_{n+1}\in A$
that fulfill (\ref{InEq.DiagProp2}) and 
$\|s_j\|^2 \leq D_{n+1}(1/\varepsilon)$.
(For general $\varepsilon_j>0,\tau>0$,
set $\varepsilon:=\min(\varepsilon_1,\dots,
\varepsilon_n)$ and decrease $\tau$ 
if $\tau>\varepsilon$.)

Base case $n=1$: Let 
$\varepsilon\ge \tau >0$ be given. We prove 
$D_2(t):=D(t;a_1,a_2)$ is a controlling function 
by finding $s_1,s_2\in A$ 
fulfilling (\ref{InEq.DiagProp2}) and $\|s_j\|^2 
\leq D_{2}(1/\varepsilon)$ for our choice of 
$D_2$.
With $x:=a_{12}$ the sequence
$y_k:= (a_1+1/k)^{-1/2} x (a_2+1/k)^{-1/2}\in A$ 
satisfies $\| y_k \| \leq 1$ and
$x:= \lim_k  x_k$  for 
$x_k:= a_1^{1/2} y_k a_2^{1/2}$,
\cf Remark \ref{rems:b-eps=d*ad}(iv).
Let  
$\delta:= \tau/ (2+ 2D_2(1/\varepsilon))$.
There is $k\in \N$ with $\| x -x_k\| <\delta$.
By assumptions on $D$, there exist
$s_1, s_2\in A$, 
with 
$\| s_j \| ^2\leq D_2(1/\varepsilon)$,
$\|  s_j^*a_j\,s_j \, -\, a_j\| 
<\varepsilon$
and 
$\| s_1^*x_k s_2\|  < \delta$.
Then 
$ \| s_1^*x s_2 \| < 
\delta + D_2(1/\varepsilon) \| x-x_k \|<\tau$,
giving (\ref{InEq.DiagProp2}). 

We proceed by induction over $n\ge 2$.
Suppose that each $n$-tuple 
$(h_1,\ldots,h_n)$ with $h_j\in \cF$ has
controlled matrix diagonalization 
with controlling functions 
${t\mapsto D_n(t;\, h_1,\ldots,h_n)}$ 
having $h_1,\ldots,h_n$ as parameters.
In particular, the functions
$t\mapsto D_2(t; a_1,a_{n+1})$, 
$t\mapsto D_n(t;a_1,\ldots,a_n)$ and 
$t\mapsto D_n(t;a_2,\ldots,a_{n+1})$
used below could be different.
We try to 
keep notations transparent by defining
$$D_2(t):=D_2(t; a_1,a_{n+1}), \ \ 
D_n(t):=
\max 
\{\,D_n(t;a_1,\ldots,a_n),\, 
D_n(t;a_2,\ldots,a_{n+1})\, \}.$$

Now let 
$\varepsilon
\ge \tau >0$ be given.
We consider the following values
\begin{equation}\label{Eq:Iteration}
\varepsilon_2:=\varepsilon/3\,, \quad
\varepsilon_1:=\varepsilon/(3D_n(3/\varepsilon))
\quad \text{and} \quad 
\varepsilon_0:= 
\varepsilon/(3D_n(1/\varepsilon_1))\,,
\end{equation}
and choose 
$\,\tau_0,\, \tau_1,\, \tau_2\, >0\,$ such that
$$
\tau_2<\tau\,,\ \ \  D_n(1/\varepsilon_2)\tau_1
<\tau\,, \ \ \
D_n(1/\varepsilon_2) D_n(1/\varepsilon_1)\tau_0
<\tau
\,.$$
Notice that $D_n(1/\varepsilon_1)=
D_n(\,D_n(1/\varepsilon_2)/\varepsilon_2\, )$
and $\varepsilon/(3D_n(3/\varepsilon))= 
\varepsilon_2/D_n(1/\varepsilon_2)$.

There are  $d_1, d_{n+1}\in A$ with 
$\| d_j \| ^2 \leq D_2(1/\varepsilon_0)$ such that
$\| d_j^*a_jd_j -a_j \| < \varepsilon_0$ for 
$j=1, n+1$ and 
$\| d_1^* a_{1,n+1} d_{n+1}\| < \tau_0$.
We can use 
Remark \ref{rems:b-eps=d*ad}(i) to modify
 $d_1$ and $d_{n+1}$  suitably, such that
that $d_j^*a_jd_j=(a_j-\mu)_+$
for some $\mu <\varepsilon_0$. 
Now consider the diagonal matrices 
$
w_1:=
\mathrm{diag}
(a_1-(a_1-\mu)_+\,,\,  0,\ldots, 0 ,\, a_{n+1} -
(a_{n+1}-\mu)_+\,)
$ 
and 
$d:= \mathrm{diag}(d_1, 1,\ldots,1,d_{n+1})$ 
in $M_{n+1}(\cM  (A))$. 
The matrix $[b_{jk}]:=w_1+ d^*[a_{jk}]d $  
is positive  in $M_{n+1}(A)$
with $b_{jj}=a_j$ and 
$b_{jk}=d_j^*a_{jk}d_k$ for $j\not=k$,
and $\| b_{1,n+1} \| < \tau_0$.

By induction hypothesis and 
Remark \ref{rems:b-eps=d*ad}(i), 
there exists a diagonal matrix 
$e=\diag (e_1, \ldots, e_n, 1)$
such that
$\| e \|^2 \leq D_n(1/\varepsilon_1)$, 
$\| e_j^*b_{j k} e_k \| <  \tau_1$ 
for 
$j\not=k \in \{ 1,\ldots, n\}$ and 
$e_j^*a_je_j=(a_j- \nu)_+$ 
for 
$j= 1,\ldots,n$ and some $\nu <\varepsilon_1$.

Consider the diagonal matrix 
$w_2:=
\diag(a_1-(a_1-\nu)_+, \ldots, a_n-(a_n-\nu)_+, 0)$
and the positive matrix 
$[c_{jk}]:= e^*[b_{jk}]e+ w_2$ 
with diagonal entries
$c_{jj}=a_j$ and $\| c_{1,k}\| < \tau_1$ for 
$k=2,\ldots,n$,  and 
$\| c_{1,k}\| <D_n(1/\varepsilon_1)\tau_0$ for 
$k=n+1$.

Apply the induction assumption to the 
lower right $n\times n$ 
sub-matrix of $[c_{j k}]$,  get a diagonal matrix 
$f:= \diag (1,f_2, \ldots, f_{n+1})$ such that
$\| f \|^2 <D_n(1/\varepsilon_2)$  and 
$$\| f_j^*a_j f_j -a_j \| <\varepsilon_2 \,,\quad 
 \quad \| f_j^*c_{jk} f_j \| < \tau_2 \quad 
 \text{for}\,\,\,\,  j\not=k \in 
\{ 2, \ldots, n+1\}\,.$$
The diagonal matrix 
$
g:= 
\mathrm{diag} 
(d_1e_1, e_2 f_2, \ldots, e_n f_n, d_{n+1} f_{n+1})
$
has norm 
$$\| g \|^2 \leq 
\max \{ D_2(1/\varepsilon_0)D_n(1/\varepsilon_1), 
D_n(1/\varepsilon_1)D_n(1/\varepsilon_2),
D_2(1/\varepsilon_0)D_n(1/\varepsilon_2) \}$$
and satisfies 
\begin{equation*}\label{xn}
\| g_j^*a_{jk} g_k\|\leq\left\{
\begin{array}{lr}
\tau_2,&\mbox{ if } j\not=k\in \{2,\dots, n+1\} \\
D_n(1/\varepsilon_2)D_n(1/\varepsilon_1)\tau_0,
&\mbox{ if }j=1,k=n+1\\
D_n(1/\varepsilon_2)\tau_1,&
\mbox{ if }j=1,k=2,\dots,n\\
\end{array}\right.
\end{equation*}
By assumption on $\tau_0,\tau_1,\tau_2$ we get 
$\| g_j^*a_{jk} g_k\|<\tau$ 
for all  $j\not=k \in \{1, \ldots, n+1\}$. Also
$$\| g_j^* a_j g_j - a_j\| <
\max\{ 
D_n(1/\varepsilon_1)\varepsilon_0 + \varepsilon_1\,,\,\, 
D_n(1/\varepsilon_2)\varepsilon_1+ \varepsilon_2\,,\,\, 
D_n(1/\varepsilon_2)\varepsilon_0 + 
\varepsilon_2\} \leq \varepsilon \,.
$$
Thus, the  $(n+1)$-tuple  
$(a_1,\ldots, a_n,a_{n+1})$ has the      
diagonalization  property with (the clearly 
non-decreasing) controlling function
$D_{n+1}(t; a_1,\ldots, a_n,a_{n+1}):=D_{n+1}(t)$ 
defined by
$$D_{n+1}(t):=
\max \{ D_2(1/\varepsilon_0)D_n(1/\varepsilon_1), 
D_n(1/\varepsilon_1)D_n(1/\varepsilon_2),
D_2(1/\varepsilon_0)D_n(1/\varepsilon_2) \}$$
with 
$(\varepsilon_0, \varepsilon_1, \varepsilon_2)$
defined from 
$\varepsilon:=1/t$ as  above  in  
(\ref{Eq:Iteration}).\qed

Recall that $\Cf_c (0,\infty]_+$ 
denotes the set of all non-negative continuous
functions $\varphi$ on $[0,\infty)$ with 
$\varphi | [0,\eta] = 0$ 
for some $\eta \in (0,\infty)$, 
such that $\lim_{t\to \infty} \varphi(t)$ exists.

In the following Lemma 
\ref{lem:from.2.diag.Prop.to.Estimate} we show 
that approximate matrix diagonalization 
for positive matrices 
in $M_2(A)$ extends in suitable cases 
to approximate 
diagonalization of (certain) selfadjoint matrices.
In particular
notice that the solvability by 
$s_1,s_2 \in \cS$ of the 
inequality (\ref{InEq.general}) 
does not
require the positivity of the
selfadjoint  $2\times 2$-matrix 
$\minimatrix{a_1}{c}{c^*}{a_2}$ anymore.

\begin{lem}
\label{lem:from.2.diag.Prop.to.Estimate}
Let $a_1,a_2\in A_+$, $\varepsilon_0 >0$
and non-empty subsets $\Co\subseteq A$, 
$\cS\subseteq A$ be given. 
Suppose that
the following
properties hold:
\begin{itemize}
\item[(i)] For every 
$\delta\in (0,\varepsilon_0)$, the pair 
$((a_1-\delta)_+, (a_2-\delta)_+)$ has 
 the matrix 
diagonalization
property with respect to $\cS$ and $\Co $
of Definition \ref{def:MatrixDiag}.
\item[(ii)]
$\varphi(a_1)c\varphi (a_2)\in \Co$ 
for each $c\in \Co$
and $\varphi \in \Cf_c (0,\infty]_+$.
\item[(iii)]
$\varphi(a_1)s, \varphi(a_2)s\in \cS$ 
for each $s\in \cS$
and $\varphi \in \Cf_c (0,\infty]_+$.
\end{itemize}
Then, for each $c\in \Co$,
$\varepsilon \in (0,\varepsilon_0)$ and 
$\tau > 0$, there exist
$s_1,s_2\in \cS$ 
that fulfill
\emph{(\ref{InEq.general})}
and 
$\|s_j\|^2 \leq 2\| a_j \| / \varepsilon$.
\end{lem}

\smallskip

\pf  
Let $\varepsilon_0 \geq \varepsilon >0$ and 
$\tau>0$, and define $\gamma:=\varepsilon/2$. 

Condition (ii) on $\Co$ implies $0\in\Co\,$.
Therefore there are $d_1,d_2\in \cS$
that satisfy the inequalities 
$\| d_j ^* (a_j-\gamma)_+ d_j - 
(a_j-\gamma)_+ \| <  \gamma$.
Then $s_j:=e_j d_j$ is in $\cS$ 
and $e_j a_j e_j= (a_j-\gamma)_+$ for 
$e_j:= \varphi(a_j)$,
where  
$\varphi(t):=( (t-\gamma)_+/t)^{1/2}\,$,
and the $\{ a_j, 0, s_j,  \varepsilon, 
\tau \}$
satisfy the general inequalities  
(\ref{InEq.general})
with $c:=0$.

Since $\gamma \varphi(t)^2\leq \varphi(t)^2t$
we get  
$s_j^*s_j\leq 
\gamma^{-1} d_j^*(a_j-\gamma)_+d_j$
and using the norm inequality,
$d_j^*(a_j-\gamma)_+d_j\leq 
\gamma +(a_j -\gamma)_+$.
Since $\| \cdot \|$ is monotone on $A_+$, 
it follows  $s_j=0$ if
$\gamma \ge  \| a_j \|$, and 
$\| s_j\|^2\leq \gamma^{-1}\| a_j \|$ 
if $\gamma <  \| a_j \|$,
so $\| s_j\|^2\leq 2\| a_j \|/\varepsilon$. 

Suppose that $\min( \| a_1\|, \| a_2 \|)=0$.
If follows that the solution $(s_1,s_2)$ 
for the above case with $c=0$ gives 
$\min( \| s_1\|, \| s_2 \|)=0$.
In particular $(s_1,s_2)$ fulfill
(\ref{InEq.general}) and 
$\| s_j\|^2\leq 2\| a_j \|/\varepsilon$ for 
each $c\in\Co$ with $\|s_1^*c s_2\|=0$.

It remains to check the case where 
$\min( \| a_1\|, \| a_2 \|)>0$. We let
$$
\alpha_1 : =\, \min( \| a_1\|, \| a_2 \|)  
\quad \text{and} \quad
\alpha_2 : =\, \max(\| a_1\|, \| a_2\|)
\,.
$$
Let  $c\in \Co$. By decreasing $\tau$ if
$\tau>\varepsilon$, we may assume 
$\tau\leq \varepsilon$.

We define functions 
$\chi, \psi, \varphi \in \Cf_c(0,\infty]_+$ by 
$\xi(t):= \min(\alpha_2, (t-\gamma)_+)^{1/2}$,
$\chi(t):= \delta^{-1}\min ((t-\delta)_+, \delta)$ 
with
$\delta:=\gamma/2$, $\psi(t):= \chi(t)t^{-1/2}$
and $\varphi$ is as above defined.
Notice that 
$\psi(t)^2(t-\gamma)_+=\varphi(t)^2$,
$\varphi(t)^2 t=(t-\gamma)_+ \,$,  
$\xi(a_j)=(a_j-\gamma)_+^{1/2}$.

Let $e_j:= \varphi (a_j)$ and
$f_j:= \psi (a_j)$.
It  follows that 
$e_j= (a_j-\gamma )_+^{1/2}f_j$ has norm
$\| e_j \| = \|\varphi (a_j)\| \leq 1$.
By assumption (iii), the elements
$f_1c f_2$ and  $e_1 c e_2$ are in
$\Co$.

Case $f_1cf_2=0\,$:  
Then  $e_1 c e_2=0$.
Since we do not know if $e_j$ is in $\cS$,
we can not define $s_j$ simply 
by $s_j:=e_j$.
But the above considered
case $c=0$ gives $d_1,d_2
\in \cS$ with
$\| d_j^*(a_j-\gamma )_+d_j 
- (a_j-\gamma )_+\| 
<\gamma $.
The $s_j:= e_j d_j\in \cS$ 
satisfy the inequalities 
(\ref{InEq.general}) and 
$\| s_j\|^2\leq \gamma^{-1} \| a_j\|$, 
where we use that 
$e_j a_j e_j= (a_j-\gamma )_+$ and above
estimates for $\|s_j\|^2$.

Case $f_1 cf_2\not=0\,$:  
We define $\rho:= \max(1,\| f_1cf_2\|)^{-1}$
and $c':=\rho \cdot f_1cf_2$
and put $\tau':= \rho \tau>0\,$.
Then $c'\in \Co$ by assumption (ii), because 
$\sqrt{\rho}\psi\in \Cf_c (0,\infty]_+$,
$\| c' \| \leq 1$ and  the matrix 
$X=[b_{ij}]\in M_2(\cM(A))$ with entries 
$b_{11}:=b_{22}:=1$
and $b_{21}^*:=b_{12}:=c' $ 
is a positive matrix
since $\|c'\|\leq 1$.
Recall that $\xi(a_j)=(a_j-\gamma )_+^{1/2}$.
Since $\xi(t)\psi(t)=\varphi(t)$ for 
$t\leq \alpha_2$, we have that 
$(a_j-\gamma )_+= \xi(a_j)f_ja_jf_j\xi(a_j)$.

Let 
$
\Delta:= 
\diag( (a_1-\gamma )^{1/2}_+,  
(a_2 -\gamma )^{1/2}_+)
\,$.
The $2\times 2$-matrix 
$[y_{ij}]=Y:=\Delta X \Delta\in M_2(A)_+$ 
has diagonal $\Delta^2$
and the upper right element of $Y$ is 
$$y_{12}=\xi(a_1)c'\xi(a_2)=
\rho (a_1-\gamma )_+^{1/2}f_1c
f_2(a_2-\gamma )_+^{1/2}
\,.$$
It is in $\Co$ by condition (ii).
By condition (i),
$(( a_1-\gamma )_+ ,  (a_2 -\gamma )_+)$ 
has the diagonalization 
property with respect to $\cS$ and  $\Co$. 
Hence for each $\mu>0$ there exist 
$d_1,d_2\in \cS$
such that the diagonal matrix 
$S =\diag(d_1,d_2)$
satisfies with respect to the norm of 
$M_2(A)$ the
inequality
$$\| S^*YS - \Delta^2 \| <\mu\,.$$
E.g.\ we can take 
$0<\mu <
\min(\tau',  \gamma, \| a_1\|, \| a_2\|)$.
Then this implies  
$$ 
\rho \cdot \| d_1^*e_1\, c \, 
e_2 d_2\|= 
\| d_1^* (a_1-\gamma)_+^{1/2} 
\,c'\, (a_2-\gamma)_+^{1/2} d_2\| 
<
\mu < \rho\cdot \tau
$$
and 
$$ \| d_j^*e_j a_j e_j d_j 
- (a_j-\gamma)_+ \| =
\| d_j^*(a_j-\gamma)_+d_j 
- (a_j-\gamma)_+ \| < 
\mu 
\,.
$$
The  $s_j:= e_j d_j$ 
($j=1,2$) are in $\cS$ by assumption (iii) and
fulfill the inequalities  
(\ref{InEq.general}). 

An upper estimate of the minimal possible 
norms of the $s_1,s_2\in A$ 
that fulfill  
the inequalities
(\ref{InEq.general}) 
can now be deduced as 
above from 
$$\gamma s_j^*s_j\leq  
d_j^*(a_j-\gamma)_+d_j\leq  
\mu+(a_j-\gamma)_+\,.
$$
It implies  that 
$
\gamma \|s_j\|^2\leq 
\mu +(\|a_j\| -\gamma)_+\leq \|a_j\|
$. \qed

Notice that the  Lemma 
\ref{lem:from.2.diag.Prop.to.Estimate} and  
Lemma \ref{lem:from.2.diag.to.matrix.diag}
together generalize 
\cite[lem.~5.6,lem.~5.7]{KirRorOinf}.
 
Combining Lemma 
\ref{lem:from.2.diag.to.matrix.diag} 
and Lemma \ref{lem:from.2.diag.Prop.to.Estimate} 
we obtain the following result applicable to 
families $\cF\subseteq  A_+$ invariant
under $\varepsilon$-cut-downs (if $\cF$ is 
not invariant
under $\varepsilon$-cut-downs one could 
enlarge $\cF$):

\begin{lem}
\label{lem:flexible.diag.imply.general}
Suppose that $\cF\subseteq  A_+$ is invariant
under $\varepsilon$-cut-downs, \ie that
for each $a\in \cF$ and 
$\varepsilon \in (0,\| a\|)$
we have
$(a-\varepsilon)_+\in \cF$.
Then the following 
conditions (i)--(iii)
on $\cF$ are equivalent:

\begin{itemize}
\item[(i)] 
Each 2-tuple $(a_1,a_2)$ 
with $a_1,a_2\in \cF$ 
has the matrix diagonalization property of 
Definition
\ref{def:MatrixDiag}.
\item[(ii)] 
For every $(a_1,a_2)\in\cF\times \cF$,
$c\in a_1^{1/2} A a_2^{1/2}$ and  
$\varepsilon  >0\, $,
there  exists $e_1,e_2\in A$ such that  
$\| e_1^*c e_2\| < \varepsilon\,$  and 
$\| a_j - e_j^*a_je_j  \| < \varepsilon $  
(for $j\in \{ 1,2 \}$).
\item[(iii)] 
For every $(a_1,a_2)\in\cF\times \cF$,
$c\in A$ and $\varepsilon >0$,
there  exists $d_1,d_2\in A$ such that  
$ \| d_1^*c d_2\| < \varepsilon$, 
$\| a_j - d_j^*a_j d_j \| <\varepsilon $ 
(for $j\in \{ 1,2\}\,$), 
and $\| d_j\|^2 \leq 2\| a_j \|/\varepsilon\,$.
\end{itemize}

In particular, the family  
$\cF$ has the controlled matrix 
diagonalization property of Definition
\ref{def5.1}, 
if and only if, each pair 
of elements in $\cF$  
has the matrix diagonalization property of 
Definition \ref{def:MatrixDiag}.
\end{lem}

\pf 
Let $(a_1,a_2)\in \cF\times \cF$.
Since 
$
\bigl((a_1-\varepsilon_1)_+,(a_2-\varepsilon_2)_+
\bigr) \in \cF\times \cF
$ 
for 
$\varepsilon_j\in [0, \|a_j\|]$,
we get that conditions (i) and (iii) are equivalent
by \cite[lem.5.6]{KirRorOinf} or 
Lemma \ref{lem:from.2.diag.Prop.to.Estimate} with
$\Co:=A$ and $\cS:=A$, which shows that
the estimate
$ \| d_j \|^2 \leq  \varepsilon^{-1}2\| a_j\|$
follows from condition (i).
Condition (iii) implies that the family $\cF$ has 
the controlled
matrix diagonalization property in the sense
of Definition \ref{def5.1} 
by Lemma \ref{lem:from.2.diag.to.matrix.diag}, 
because
for $a_1,a_2\in \cF$ we can define 
a control  $t\mapsto D_2(t;a_1,a_2)$ 
by 
$D_2(t; a_1,a_2):=
2\max(\| a_1\|, \| a_2\|) \cdot t$.
Notice that 
$2\| a_j \|/\varepsilon \leq 
D_2(1/\varepsilon;a_1,a_2)$ for $j\in \{ 1,2\}$.
The controlled matrix diagonalization property 
of the family $\cF$ with respect to $\cS:=\Co:=A$
follows from
Lemma \ref{lem:from.2.diag.to.matrix.diag}
if $\cF\times \cF$ satisfies condition (iii).
Clearly condition (iii) implies condition (ii)
on $\cF$.
It remains to show that condition (ii) implies 
condition (i) for $\cF\times \cF$.

(ii)$\Rightarrow$(i): 
By Remark \ref{rems:b-eps=d*ad}(iv), 
the condition (ii) on the 
pairs $(a_1,a_2)\in \cF$
is formally weaker than the matrix-decomposition
property. 
But the fact that condition (ii)
applies also to
$\bigl(\,(a_1-\delta)_+,(a_2-\delta)_+\bigr)
\in \cF\times \cF$
-- in place of $(a_1,a_2)$ -- 
allows to proof
that $(a_1,a_2)$ has the 
matrix diagonalization property, and that
\cite[lem.~5.6]{KirRorOinf} 
(or alternatively  
Lemma \ref{lem:from.2.diag.Prop.to.Estimate}
with $\Co:=(a_1-\delta)_+^{1/2}A(a_2-\delta)_+^{1/2}$
and $\cS:=A$ for each $\delta>0$)
applies
and shows the estimate given in
condition (iii), which is (i) and the estimate.

More generally,
let $a_1,a_2\in \cF\subset A_+$.
Then following observation on 
general  $a_1,a_2 \in A_+$  
applies to $(a_1,a_2)$ and shows
that $(a_1,a_2)$ has the matrix diagonalization
property of Definition \ref{def:MatrixDiag}.
It follows then that each of the conditions (i)--(iii)
on $\cF\times \cF$ are equivalent.

\emph{Suppose that 
$a_1,a_2\in A_+$ are non-zero
and have the property
that,   
for each 
$\delta,\tau \in (0,\min\{ \|a_1\|,\|a_2\|\})$ 
and each $x\in A$,
there exist $e_1,e_2\in A$,
depending on $(x,\delta,\tau)$,
with 
%$$
\begin{equation}\label{InEq.with.cutdown}
\|e_j^*(a_j-\delta)_+ e_j -(a_j-\delta)_+\| 
<\tau
\quad \text{and} \quad
\| e_1^*(a_1-\delta)_+^{1/2}x
(a_2-\delta)_+^{1/2} e_2\| <\tau\,.
\end{equation} 
%$$
Then there exist for each $\varepsilon>0$
and $c\in A$
elements $d_1,d_2\in A$
with $\| d_j^*a_j d_j -a_j\|<\varepsilon$, 
$\,j\in \{1,2\}$,
and $\|d_1^*c d_2\| < \varepsilon$.
}

\smallskip

To verify this Observation,
let $a_1,a_2\in A_+$, $c\in A$ and
$\varepsilon>0$ given. 
Without loss of generality  we can suppose
that $a_1,a_2$ and $c$ are non-zero
and that $2\varepsilon\leq \min(\|a_1\|,\|a_2\|)$.

Let
$\tau:=\gamma:=\varepsilon/2$,
and $\delta\in (0,\gamma)$. Then
$$\| (a_j-\delta)_+ - a_j \|\leq \delta < \gamma$$
and there exists a bounded continuous real function 
$g_\delta\in \mathrm{C}_c(0,\infty]_+$
with $g_\delta(\xi)=0$ on $[0,\delta/2]$
and $g_\delta(\xi)=\xi^{-1/2}$
for $\xi \ge \delta$.
Let $h_\delta(\xi):= \xi^{-1/2}(\xi-\delta)_+^{1/2}$
for $\xi\in (0,\infty)$ and $h_\delta(0):=0$.

Notice 
that $\xi ^{1/2} h_\delta (\xi)=(\xi-\delta)_+^{1/2}$
and $g_\delta(\xi) (\xi-\delta)_+^{1/2}= h_\delta(\xi)$.
Let $x:= h_\delta(a_1)c h_\delta(a_2)$.
By assumptions of the Observation, 
there exist $e_1,e_2\in A$
with 
$$
\| e_j^*(a_j-\delta)_+e_j -(a_j-\delta)_+\|
< \tau \quad \text{and} \quad
\| e_1^*x e_2\| < \tau
\,. $$
Now let $d_j:= h_\delta(a_j)e_j$.
Then $e_1^*x e_2=d_1^* c d_2$
and $d_j^*a_j d_j = e_j^*(a_j-\delta)_+e_j$.
Thus, 
$\| d_j^*a_j d_j -a_j \| \leq \delta +\gamma< \varepsilon$
and $\| d_1^*cd_2\| \leq \tau <\varepsilon$, as desired.
\qed

In the following two Lemmata 
\ref{lem:passage.to.span.C} and 
\ref{lem:cS.control}
we consider a
globalization of
Lemma \ref{lem:from.2.diag.Prop.to.Estimate}
to the case of families 
$\cF\subseteq A_+$ and the case where
$\cS$ is moreover  
a multiplicative sub-semigroup of $A$ 
satisfying $s_2^*\Co s_1\subseteq  \Co$ 
for all $s_1,s_2\in \cS$. We introduce 
first a bit notation. Let 
$D\colon (0,\infty)\to [1,\infty)$ be any
function and $n$ any positive integer.
For each $t\in (0,\infty)$
let $\cY_{t,n}$, or simply $\cY_t$, 
denote set\footnote{\,The set $\cY_t$ is
nonempty for every $t>0$ by Proof of 
Lemma \ref{lem:passage.to.span.C}.} 
consisting of all $n$-tuples 
$(\varepsilon_1,\ldots,\varepsilon_n)
\in (0, 1/t)^n$
that satisfy the inequality
\begin{equation}\label{InEq.Yt.condition}
\varepsilon_n \, +\,\,
\sum _{k=1}^{n-1} \,\, \varepsilon_k \cdot 
D(1/\varepsilon_{k+1})
\cdot \ldots \cdot D(1/\varepsilon_n)
\, \leq \, 1/t\,.
\end{equation}
Moreover, let  $t\mapsto E_n(t)\in [1,\infty)$ 
denote the function defined by 
$$ 
E_n(t):= \inf\{\, D(1/\varepsilon_1)\cdot 
\ldots \cdot 
D(1/\varepsilon_n)\,;\,\,  
(\varepsilon_1\,,\,\ldots\,,\, 
\varepsilon_n)\in \cY_t\}
\,.$$ 

\begin{lem}\label{lem:passage.to.span.C}
Let $a_1,a_2\in A_+$, $\varepsilon_0\in (0,\infty]$
and non-empty subsets $\cS\subseteq A$, 
$\Co\subseteq A$
be given. 
Suppose that
the following
properties hold:
\begin{itemize}
\item[(i)] The set $\,\cS$ is a
multiplicative sub-semigroup of $A$ that 
satisfies $s_2^*\Co s_1\subseteq  \Co$ 
for all $s_1,s_2\in \cS$.
\item[(ii)] There exists a non-decreasing
function
$t\in (0,\infty) \mapsto 
D(t):=D(t ; a_1,a_2)\in [1,\infty)$
such that
for each $c\in \Co$, each 
$\tau,\varepsilon \in (0,\varepsilon_0)$ 
with $\varepsilon  \geq  \tau$
there exist 
$s_1,s_2\in \cS$
that fulfill \emph{(\ref{InEq.general})}
and 
$\|s_j\|^2 \leq D(1/\varepsilon)$.
\end{itemize}	
	
Then one can find, for each finite subset 
$X\subseteq \Co$
and $\varepsilon\ge \tau >0$,
elements $s_1, s_2\in \cS$ 
that satisfy \emph{(\ref{InEq.general})}
for \emph{every} $c\in X$.
Moreover, if 
$\varepsilon\in (0, \varepsilon_0)$ and 
$n:= |X|$, then we can ensure 
$\| s_j \|^2\leq E_n (1/\varepsilon)$.

\smallskip
For each 
$c\in A$  in the 
linear span of at most $n$ elements 
$c_1,\ldots, c_n\in \Co$
and each $\varepsilon\ge \tau >0$,
there exist
$s_1,s_2\in \cS$ 
that fulfill
\emph{(\ref{InEq.general})}.
Moreover, 
if $\varepsilon\in (0,\varepsilon_0)$, then 
we can ensure 
$\| s_j \|^2\leq E_n (1/\varepsilon)$.

\smallskip
If $D(t) \leq \gamma$ for a constant $\gamma$
then $E_n (t) \leq \gamma^n $ for all $n\in \N$, 
and if $D(t) \leq \gamma \cdot t$ 
then an upper estimate 
for $E_n$  is given by 
$E_n(t)\leq (n t \gamma)^{(2^n-1)}$.

If $D(t)=1$ then, 
for each $c$ in the closure of the linear span of 
of $\Co$, $\varepsilon\in (0,\varepsilon_0)$, and $\tau>0$, 
there are contractions 
$s_1,s_2\in \cS$ that satisfy the 
inequalities \emph{(\ref{InEq.general})}. 
\end{lem}

\pf
If $(a_1,a_2)$, $\Co$, $\cS$, 
$\varepsilon_0$ 
and $t\mapsto D(t)$
are given,  
then we can define for any positive $t \in \R$
numbers $\nu_0,\nu_1,\ldots,\nu_{n-1}$
by induction for $k=0,1,\ldots, n-1$ as follows:
Let $\nu _0 := 2t$ and 
$\nu_{k+1}:= 2D(\nu_k)\nu_k$.
The $n$-tuple 
$(\varepsilon_1,\ldots,\varepsilon_n)$
with 
$\varepsilon_k:= \nu_{n-k}^{-1}$ ($k=1,\ldots n$)
satisfies the inequality 
(\ref{InEq.Yt.condition})
with ``$<$'' in place of ``$\leq$''.
Thus $\cY_t$ is non-empty.
An \emph{alternative construction} is given by 
$\nu_0:=n t$ and 
$\nu_{k+1}:= D(\nu_k)\nu_k$.
Then the $\varepsilon_k:=\nu_{n-k}^{-1}$
satisfy  (\ref{InEq.Yt.condition}) with ``$=$''
in place of ``$\leq$''. 
We use the latter to find
bounds for $E_n(t)$.

Given 
$X= \{ x_1,\ldots, x_n\} \subseteq  \Co$ 
and
$\varepsilon \geq \tau >0$.
If $\varepsilon\ge \varepsilon_0$ then
we can decrease  $\varepsilon$ 
and $\tau$ below $\varepsilon_0$
in the following considerations and
find the $s_1,s_2$ 
satisfying \eqref{InEq.general} for these
smaller $\varepsilon$ and $\tau$.
We only need  to prove 
the norm estimate  
$\|s_j\|^{2} \leq  E_n(1/\varepsilon)$
for $\varepsilon \in (0,\varepsilon_0)$.
Therefore we can assume below
always that $\varepsilon<\varepsilon_0$.

With $t:=1/\varepsilon$,
$\cY_t$ is non-empty
by the computation above.
Let $(\varepsilon_1,\ldots,\varepsilon_n)$
an \emph{arbitrary} element of $\cY_t$
and define $\tau_k$ from $\tau$ and the 
$\varepsilon_k$
by
$$\tau_k:= \min\{
\tau / D( 1/ \varepsilon_{k+1} ) 
\cdot \ldots \cdot  
D(1/\varepsilon_n), \varepsilon_k\}\,.$$
By assumptions applied to 
$\tau_k,\varepsilon_k\in (0,\varepsilon_0)$, 
we can find 
elements
$s_1^{(k)},s_2^{(k)}\in \cS$  
with norms $\| s_j^{(k)} \| \leq 
D(1/\varepsilon_k)^{1/2}$ 
that satisfy the inequalities 
(for $j=1,2$, $k=1,\dots,n$)
\begin{equation}\label{sk.estim}
\| (s_j^{(k)}) ^* a_j s_j^{(k)} - a_j \| < 
\varepsilon_k 
\quad \text{and} \quad
\| (s_1^{(k)})^*\, c_k \,s_2^{(k)}\| < \tau_k
\,,\end{equation}
where we let $c_1:= x_1$ and
$c_{k+1}:= 
\bigl( s_1^{(1)} \cdots s_1^{(k)} \bigr)^* x_{k+1} 
\bigl( s_2^{(1)} \cdots s_2^{(k)} \bigr) \in \Co
\,$.

The  $s_j := s_j^{(1)} \cdots s_j^{(n)}$ 
($j=1,2$) satisfy
\begin{equation}\label{InEq.dj.bound}
\| s_j\|^2\leq D(1/\varepsilon_1)\cdot \ldots 
\cdot D(1/\varepsilon_n)
\end{equation}
and
\begin{equation*}
\| s_1^* x_k  s_2 \| < 
\tau_k \cdot  D( 1/ \varepsilon_{k+1} ) 
\cdot \ldots \cdot  
D(1/\varepsilon_n) \leq \tau
\,.
\end{equation*}
Stepwise application of the triangle inequality 
and (\ref{sk.estim})
shows that
$$
\| s_j^*a_j s_j - a_j\| <  
\varepsilon_n + \sum_{k=1}^{n-1} \varepsilon_k 
D(1/\varepsilon_{k+1})\,\cdot\, \ldots\, \cdot\, 
D(1/\varepsilon_n) 
\,.$$
Since $(\varepsilon_1,\ldots,\varepsilon_n)
\in \cY_{1/\varepsilon}$ we get
$\| s_j^*a_j s_j - a_j\| < 
\varepsilon$, ensuring (\ref{InEq.general}) 
for each $c=x_k$.
Since $(\varepsilon_1,\ldots,\varepsilon_n)
\in \cY_{1/\varepsilon}$ was
arbitrary, $\| s_j\|^2\leq E_n(1/\varepsilon)$.

The passage to the $c$ in the linear span of 
finite subsets 
$X:= \{ x_1,\ldots, x_n\}\subseteq  \Co$
is a matter of finding a solution 
$s_1,s_2$ of the inequalities (\ref{InEq.general}) 
for all $c\in X$ with appropriate choice of $\tau$,
(for $c=\sum_{i=1}^n \alpha_i x_i$ 
selecting $s_j$ such that $\| s_1^* x_k  s_2 \|
\leq\tau/(1+ \sum_{i=1}^n |\alpha_i|)$ will do
the job). The results on the norm estimates 
on $s_1,s_2$ remains
unchanged.

The estimate  $E_n(t)\leq \gamma^n$ for 
constant $D(t)=\gamma$ follows easily 
by the definition of $E_n(t)$,
because we know that $\cY_t$
is non-empty.
If $D(t)=\gamma\cdot t$, $t\ge 1$, 
then take the $\nu_k$ as in the above 
alternative construction.
It follows $(\frac{1}{\nu_{n-1}},\dots,
\frac{1}{\nu_{0}})\in \cY_t$,
hence $E_n(t)\leq D(\nu_0)\cdot\, \ldots\, 
\cdot D(\nu_{n-1}) 
= 
\nu_0 \cdot\, \ldots\, \cdot \nu_{n-1}
\gamma^n$, and 
$\nu_k = (nt\gamma)^{2^k}/\gamma$. 
Thus $E_n(t)\leq (nt\gamma)^{2^n-1}$.

If $D(t)=1$ then $E_n(t)=1$ for all $n\in \N$.
Hence we can decrease $\varepsilon$ without
enlarging $E(1/\varepsilon)$.
Consequently, for 
$c$ in the closed linear span 
of $\Co$ and \emph{any}
$\varepsilon\in (0,\varepsilon_0)$,
$\tau>0$, we can
find $c'\in \mathrm{span}(\Co)$ with 
$\| c'-c \| <\tau/2$  and 
contractions $s_1,s_2\in \cS$ that
satisfy the inequalities (\ref{InEq.general})
with $(c',\tau/2)$ in place of $(c,\tau)$.
Then this $s_1,s_2$ also satisfy 
(\ref{InEq.general})
with the given $c$ and $\tau$.
\qed

\begin{lem}\label{lem:cS.control}
Let $\varepsilon_0>0$
and non-empty subsets $\cF\subseteq A_+$,
$\Co\subseteq A$
be given, and let 
$\cS \subseteq A$ be a 
(multiplicative) sub-semigroup of $A$ 
that satisfies $s_2^*\Co s_1\subseteq  
\Co$ for all $s_1,s_2\in \cS$.
Suppose that
the following
properties holds for every $a_1,a_2\in \cF$:
\begin{itemize}
\item[(i)] 
For every $\delta\in (0, \varepsilon_0)$, the pair has
$((a_1-\delta)_+, (a_2-\delta)_+)$  
 the matrix 
diagonalization
property with respect to $\cS$ and $\Co $
 of Definition 
\ref{def:MatrixDiag}.
\item[(ii)]
$\varphi(a_1)c\varphi (a_2)\in \Co$ 
for each $c\in \Co$
and $\varphi \in \Cf_c (0,\infty]_+$.
\item[(iii)]
$\varphi(a_1)s, \varphi(a_2)s\in \cS$ 
for each $s\in \cS$
and $\varphi \in \Cf_c (0,\infty]_+$.
\end{itemize}
Then, for every 
$c\in \overline{\mathrm{span}(\Co )}$,
$a_1,a_2\in \cF$,
$\varepsilon \in (0,\varepsilon_0/2)$, and
$\tau>0$, 
there exists 
$s_1,s_2\in \cS$
that fulfill
\emph{(\ref{InEq.general})}
and 
$\|s_j\|^2 \leq 2\| a_j \|/\varepsilon$.
\end{lem}

\pf  Take any $a_1,a_2\in \cF$ and 
$\delta \in (0,\varepsilon_0/2)$. Due to 
property (i) the system 
$((a_j-\delta)_+, \varepsilon_0/2, \Co, \cS)$ 
fulfills the conditions of Lemma 
\ref{lem:from.2.diag.Prop.to.Estimate}. 
Hence, by Lemma 
\ref{lem:from.2.diag.Prop.to.Estimate}, 
for every $c\in \Co$, 
$\varepsilon \in (0, \varepsilon_0/2)$ and 
$\tau > 0$, 
there exist $s_1,s_2\in \cS$ with norms 
$\|s_j\|^2 \leq  2\|(a_j-\delta)_+\| /\varepsilon\,$ 
that satisfy the inequalities 
$\| \, s_j^*(a_j-\delta)_+s_j-(a_j-\delta)_+\,
\|<\varepsilon$, 
and $\| \, s_1^*c s_2\,\| < \tau$. Define 
$\varepsilon_1:=\varepsilon_0/2$, and 
$D_1(t):=\max(1, t\cdot 2\max \|(a_j-\delta)_+\|)$. 
The system $((a_j-\delta)_+, \Co, \cS)$ fulfills 
conditions of Lemma \ref{lem:passage.to.span.C} 
using ($\varepsilon_1,D_1$) in place of 
($\varepsilon_0,D$). 

Hence, by Lemma 
\ref{lem:passage.to.span.C}, for every 
$c\in \mathrm{span}(\Co)$ and 
 $\varepsilon\ge \tau >0$, 
there exist $s_1,s_2\in \cS$ with 
$\| \, s_j^*(a_j-\delta)_+s_j-
(a_j-\delta)_+\,\|  \, <\varepsilon$, and 
$\| \, s_1^*c s_2\,\|  \, <\tau$. We 
obtain that the system 
$(a_j, \varepsilon_0/2, \mathrm{span}(\Co), \cS)$ 
fulfills the conditions of Lemma 
\ref{lem:from.2.diag.Prop.to.Estimate}. We can 
now apply Lemma 
\ref{lem:from.2.diag.Prop.to.Estimate} 
on arbitrary $a_1,a_2\in \cF$.
It follows that for $a_1,a_2\in \cF$, 
$\varepsilon \in (0,\varepsilon_0/2)$,
$\tau>0$, 
and $c\in \mathrm{span} (\Co )$ there exist
$s_1,s_2\in \cS$ with 
$\|s_j\|^2\leq 2\| a_j \|/\varepsilon$
that satisfy the inequalities 
(\ref{InEq.general}).

If $a_1,a_2\in \cF$,
$\varepsilon \in (0,\varepsilon_0/2)$, and
$\tau>0$ are given
and
if $c=\lim _n  c_n$  with 
$c_n\in \mathrm{span}(\Co)$,
then we find $k\in \N$  with 
$\|c_k -c \| <\delta$, where 
$\delta:=
(\varepsilon \cdot \tau)/
(2 + 4\sqrt{\| a_1\|\cdot \| a_2\|})$. 
We find $s_1,s_2\in \cS$ with 
$\| s_j\|^2\leq 2\| a_j \| /\varepsilon$
that satisfy the inequalities 
(\ref{InEq.general}) with 
$(c_k,\tau/2)$ in place of 
$(c,\tau)$.
Then 
$\| s_1^*c s_2\| < 
\tau/2  +  \delta \|s_1\| \|s_2\| 
\leq  \tau$.
Hence, for given $a_1,a_2\in \cF$, 
$c\in \overline{\mathrm{span}(\Co)}$,
$\varepsilon \in (0,\varepsilon_0/2)$, and
$\tau>0$, there exist 
$s_1,s_2\in \cS$ that satisfy 
(\ref{InEq.general})
and have norms 
$\| s_j \|^2  \leq 2\| a_j \| /\varepsilon$.
\qed

\section{Pure infiniteness
of tensor products}\label{sec6}
The following 
Lemma \ref{lem:B.ot.C.matrix.diag}  
considers a subset $\cF\subseteq  A_+$
that is not invariant under 
$\varepsilon$-cut-downs.
Therefore we use the definition of 
s.p.i.\ \cst-algebras that predicts that 
inequalities (\ref{InEq.spi.final}) can be solved
by \emph{contractions} $s,t$, see 
Remark \ref{rem:dn.contraction}.
One could also work with weaker estimates 
for the $\| d_j \|$  that we can derive
with our methods here, but that would require to
prove first a more complicate version of the
local Lemma 
\ref{lem:from.2.diag.Prop.to.Estimate}
and then to update its  
globalization in Lemma \ref{lem:cS.control}.
The key observation for such a 
generalization had to be started with
comparison of $\varepsilon$-cut-downs with
elements that actually exist in $\cF$ and
how these multiply
$\cS$ from the left and $\Co$ form both sides.
The reader can find such a generalization 
for contractions 
$b\in B_+$, $c\in C_+$ and 
$0\leq \varepsilon \leq 1$
\eg with help of
the inequalities 
$((b\ot c)  - 3\varepsilon)_+ \leq 
(b-\varepsilon)_+\ot (c-\varepsilon)_+\leq 
((b\ot c) -\varepsilon^2)_+\,$.
\begin{lem}\label{lem:B.ot.C.matrix.diag}
Suppose that at least one of the 
\cst-algebras $B$ or $C$ is s.p.i.\  
Then the family
$\cF:=\{ b\otimes c\,;\,\, b\in 
B_+\,,\,c\in C_+\,\}$
has the matrix diagonalization property
in $B\otimes ^{\alpha} C$ for each \cst-norm 
$\| \,\cdot \,\|_\alpha$ 
on the algebraic tensor product
$B\odot C$.
\end{lem}
 
\pf We consider the case where 
$B$ is s.p.i.
The case of s.p.i.\  $C$ is similar.

If $b_1,b_2\in B_+$, $f\in B$,
$c_1,c_2\in C_+$,  $g\in  C$,
and $\varepsilon\ge \tau >0$ are given, we
let 
$\delta:=\,\tau/ \bigl(1+
\max\{\, \|b_1\| +\| c_1\|,\,\|b_2\| +
\| c_2\|,\,\| g\|\,\} 
\bigr)$.
By Remark \ref{rem:dn.contraction}
there exists contractions $d_1,d_2\in B$
such that $\| d_j^*b_jd_j -b_j\| <\delta$
and $\| d_1^*fd_2\| <\delta$.
Since the set of positive contractions in $C$
contains an approximative unit for $C$
(\cf  \cite[thm.~1.4.2]{Ped.book}),
there exists a contraction $e\in C_+$
with $\| e c_j e - c_j\| <\delta$.
The tensors $s_j:= d_j\ot e$ satisfy
$\| s_j^*(b_j\ot c_j)s_j -b_j\ot c_j \| 
<\varepsilon$ and 
$\| s_1^*(f\ot g) s_2\| < \tau$.

It follows that 
the pair
$(a_1, a_2)$  with $a_j:=b_j\ot c_j\in \cF$, 
the subset 
$\Co:=\{ b\otimes c\,;\,\, b\in B,\, c\in C\,\}$
of $A:=B\ot ^{\alpha} C$
and the multiplicative sub-semigroup
$\cS:= 
\{  s \otimes e \,;\,\,  
s\in B,\, e\in C,\, \| s\|\leq 1,\, \| e \|
\leq 1 \,\}$
of the algebraic tensor product
$B\odot C \subseteq  B\ot ^{\alpha} C$
satisfy the
assumptions of 
Lemma \ref{lem:passage.to.span.C} with
$\varepsilon_0:=+\infty$.
Since $\cS$ consists of
contractions, the corresponding estimating
function is $D(t)=1$.
The closed linear span of $\Co$ is dense
in $B\otimes ^{\alpha} C$. Now Lemma
\ref{lem:passage.to.span.C} gives 
that for each $a_1,a_2\in \cF$,
$c\in B\otimes ^{\alpha} C$, $\varepsilon>0$
and $\tau>0$ there exist contractions 
$s_1,s_2\in \cS$
that satisfy the inequalities 
(\ref{InEq.general}).
Therefore, Lemma 
\ref{lem:from.2.diag.to.matrix.diag}
applies to $\cF$ and we obtain that
$\cF$ has the diagonalization property
in $A$.
\qed

\begin{lem}\label{lem:a.otimes.b.filling}
It $B$ and $C$ are \cst-algebras where
$B$ or $C$ is exact, then the subset
$\cF = 
\{ b\otimes c \,;\,\, b\in B_+, c\in C_+ \}$
of  $(B \ot ^{\min} C)_+$ is a filling family
for $B \ot ^{\min} C$.
\end{lem}

\pf
Suppose that one of the algebras $B$ or 
$C$ is an exact \cst-algebra,
that $D$ is a hereditary \cst-subalgebra of 
$B\ot ^{\min} C$ and
that $I$ a primitive  ideal of 
$B \ot ^{\min} C$ with $D\not\subseteq I$.
Then  \cite[prop.~2.16(iii)]{BlanKir2},  
\cite[prop.~2.17(ii)]{BlanKir2}
and \cite[lem.~2.18]{BlanKir2} 
together show that
there exist non-zero $g\in B_+$, $h\in C_+$,  
$t\in B\ot^{\min} C$ and pure states 
$\varphi$ on $B$ and $\psi$ on $C$ such that
$(\varphi \ot \psi) (I)=\{ 0\}\,$,
$t t^*\in D$, $t^*t = g\ot h$, 
$\varphi (g) = \| g\| = 1$ and 
$\psi (h) = \| h \| = 1$.

Thus, the subset  
$\cF = 
\{ b\otimes c \,;\,\, b\in B_+, c\in C_+ \}\subseteq  
(B \ot ^{\min} C)_+$
satisfies the property (ii) of 
Lemma \ref{lem:filling.family} for all primitive 
ideals $I$ of $B\ot ^{\min} C$ with 
$D\not\subseteq I$.
By Definition \ref{def:filling.mdiag.family}
the set $\cF$ is a filling family for 
$B\ot ^{\min} C$.
\qed

\textbf{Proof of 
Theorem 
\ref{thm:A.ot.B.spi.if.A.spi.B.exact}}:\,
The subset 
$\cF:= 
\{ a\ot b\,;\,\, a\in A_+, b\in B_+\} \subseteq 
 (A\ot^{\min} B)_+$
is a filling family for $A\ot^{\min} B$ 
by Lemma 
\ref{lem:a.otimes.b.filling}, 
and 
has the matrix diagonalization property in 
$A\ot^{\min} B$
by Lemma \ref{lem:B.ot.C.matrix.diag}. 
It follows that $A\ot^{\min} B$ is
s.p.i.\  by Theorem \ref{prop:local-spi}.
\qed 

\begin{ex}\label{ex:A.ot.max.B.not.spi}
The statement of 
Theorem \ref{thm:A.ot.B.spi.if.A.spi.B.exact}
does not hold for the 
maximal tensor product of \cst-algebras:
Let 
$A:=R\ot^{\max} C^*_\lambda (F_2)$, 
$B:= C^*_\lambda (F_2)$,
where $R$ denotes the 
stably infinite simple unital nuclear 
\cst-algebra 
with finite unit element 
constructed by M.~R{\o}rdam 
 \cite{Rordam.Example},
and $F_2$ is the free group on 
two generators.
The algebras  $A$ and $B$ are exact, and 
$A$ is s.p.i.\ by 
\cite[cor.~3.11]{BlanKir2}.
The maximal \cst-tensor product 
$A\ot^{\max} B$ is even not locally purely
infinite (\cf \cite{BlanKir2} for a definition), 
because $R\ot^{\max} \cK$ 
is an ideal of a quotient of $A\ot^{\max} B\,$. 
This follows
from the fact that
the \cst-algebra generated by the ``two-sided''
regular representation 
$(g,h)\mapsto \lambda(g)\rho(h)$
of  $F_2\times F_2$ on $\ell_2(F_2)$
contains the compact operators 
in its closed linear span, 
\cf 
\cite{Akemann.Ostrand.1975.Cst.F2xF2}.
\end{ex}

\section{Endomorphism crossed product}
\label{sec:endo.cross}
Let $\varphi\colon A\to A$  
be a *-endomorphism of a
\cst-algebra $A$ that is 
\emph{not necessarily injective}.
We let 
$A_\infty:=\ell_\infty(A)/c_0(A)$ and
denote by $(A_e,\Z,\sigma)$ the 
canonical 
\cst-dynamical system associated with 
$\varphi$.

More precisely,  we consider the  inductive
limit 
$(A_e,\, \varphi_m\colon A\to A_e)$ 
 in the \emph{category} of 
 \cst-algebras of the sequence
$\xymatrix{A \ar[r]^\varphi & A \ar[r]^\varphi & 
A \ar[r]^\varphi & \cdots }\,$.
The natural realization  of 
$A_e$ is the closure of
the set 
$\bigcup _m \varphi_m(A)$ in 
$A_\infty$,
where 
$\varphi_1(a):= 
(a, \varphi(a), \varphi^2(a),\ldots)+c_0(A)$
and 
where 
$\varphi_n(a):= 
S_+^{n-1}(a,\varphi(a),\varphi^2(a),\ldots) 
+c_0(A)$
for the forward shift    
$S_+(a_1,a_2,\ldots):=(0,a_1,a_2,\ldots)$
on $\ell_\infty(A)$.

The backward shift 
$\sigma\colon (a_1,a_2,\ldots )+c_0(A)
\mapsto (a_2,a_3,\ldots )+c_0(A)$ 
is an
automorphism of $A_\infty=\ell_\infty(A)/c_0(A)$
that  induces an automorphism 
$\sigma|A_e \colon A_e \to A_e$
because 
$\varphi_n=\sigma \circ \varphi_{n+1}$ 
and $A_e$ is the closure of 
the union of
$\varphi_1(A)\subseteq 
\varphi_2(A)\subseteq \dots$. 
We denote the restriction $\sigma|A_e$
of $\sigma$ to $A_e$ simply again by $\sigma$. 
The corresponding
$\Z$-action given by $n\mapsto \sigma^n$ 
will be also denoted by 
$\sigma$
and is usually called 
\emph{the
action of the integers 
$\Z$ on $A_e$ corresponding to $\varphi$.}

The *-homomorphisms 
$\sigma, \varphi_n$ and 
$\varphi$ satisfy the equations
$$
\sigma \circ \varphi _n = 
\varphi^\infty \circ \varphi_n
\quad, \ \ \ \varphi_{n} = 
\varphi_{\ell} \circ \varphi^{\ell-n}
=
\sigma^{\ell-n}\circ \varphi_\ell \quad 
\text{and} \quad 
\sigma^k \circ \varphi_n=
\varphi_n \circ \varphi^k
\,,$$
where $1\leq n \leq \ell$, and where 
$\varphi^\infty ((a_1,a_2,\ldots)+c_0(A)):=
(\varphi(a_1),\varphi(a_2),\ldots)+c_0(A)\,$. 

Another explanation for these formulas 
can be seen from the formulas given by Cuntz 
in \cite[p.~101]{Cuntz:Symposia38.1982}
for the restriction of 
$\varphi^\infty $ to $A_e$ by the
commuting diagram:

$$
\xymatrix{
A \ar[r]^\varphi \ar[d]^\varphi & 
A \ar[r]^\varphi \ar[d]^\varphi & 
A \ar[r]^\varphi \ar[d]^\varphi & \cdots \\
A \ar[r]^\varphi  & 
A \ar[r]^\varphi & A \ar[r]^\varphi  & 
\cdots 
}
\,.$$

Recall that $\cM(A)$ 
denotes the (two-sided)
multiplier algebra of $A$.
Consider the following
non-degeneracy property (ND)
and corner property (CP):
\begin{itemize}
\item[(ND)] $A\not=I_\varphi$ for
$I_\varphi :=  
\overline{\bigcup_n  (\varphi^{n})^{-1} (0)}$.
\item[(CP)] 
The hereditary \cst-subalgebra
$\varphi(A)A\varphi(A)$ of 
$A$ is a corner of $A$.
Equivalently: $\varphi$ extends to a
strictly continuous *-homomorphism 
$\cM (\varphi)\colon \cM(A) \to    \cM(A)$.
\end{itemize}

\begin{lem}\label{lem:corner.image}
Let $\varphi\colon A\to \cM(B)$ a 
*-homomorphism and consider the 
hereditary \cst-subalgebra 
$D:=\varphi(A)B\varphi(A)$ of $B$.
The following are equivalent:
\begin{itemize}
\item[(i)] $\varphi$ extends to a strictly continuous
*-homomorphism 
$\cM(\varphi)\colon \cM(A)\to \cM(B)$.
\item[(ii)] $D$ 
is a corner of $B$, i.e., there is a projection 
$p\in \cM(B)$ with
$pBp=D$.
\item[(iii)] The (two-sided) annihilator 
$\mathrm{Ann}(D,B):=\,
\{ b\in B\,;\,\, bD=\{ 0\}=Db\, \}
$
of $D$ in $B$ has the property that
the \cst-subalgebra 
$D+ \mathrm{Ann}(D,B)$
of $B$ contains an approximate unit of $B$.
\end{itemize}
\end{lem}

\pf
It is well know and easy to see that if 
$C\subseteq \cM(B)$ is a \cst-subalgebra 
such that $CBC=B$, then there is a 
unique unital \cst-morphism $\psi$ from 
$\cM(C)$ into $\cM(B)$ that extends the 
inclusion map $C\hookrightarrow \cM(B)$ 
(i.e., $\psi(c)=c$ for $c\in C$) 
and is strictly continuous with respect 
to the strict topologies on 
$\cM(C)$ and $\cM(B)$.

(i)$\Rightarrow$(ii,iii):\, 
Since $\cM(\varphi)$ is strictly continuous, 
we get for $p:=\cM(\varphi)(1_{\cM(A)})$ that 
$D=pBp$, and $pBp+(1-p)B(1-p)$ contains 
an approximate unit, e.g. 
$pe_\tau p+(1-p)e_\tau(1-p)$ if $e_\tau$ 
is an approximately central with $p$ 
commuting with the approximate unit of $B$.

(iii)$\Rightarrow$(ii):\, Suppose that 
$E:=D+\mathrm{Ann}(D,B)$ contains an 
approximate unit of $B$, then $EBE=B$. 
Thus there exist $\psi\colon \cM(E)\to \cM(B)$ 
unital and strictly continuous, with $\psi(c)=c$ 
for $c\in E$. Since $D\cdot \mathrm{Ann}(D,B)=\{0\}$ 
we get that $E=D+ \mathrm{Ann}(D,B)$ 
is naturally isomorphic to 
$D\oplus \mathrm{Ann}(D,B)$. 
Thus 
$\cM(E)\cong \cM(D)\oplus\cM(\mathrm{Ann}(D,B))$ 
and $p:=\psi(1_{\cM(A)})$ 
is a projection in 
$B$ with $p d p=\psi(d)$ for $d\in D$.
If $b=p b p$, then 
$b=\lim d_\tau b d_\tau$ for an 
approximate unit 
$(d_\tau)$ of $D$. Thus $pBp=D$.

(ii)$\Rightarrow$(i):\, 
If $p\in \cM(B)$ is a projection satisfying 
$pBp=D$ then $\cM(D)\cong p\cM(B)p$ 
by a natural isomorphism that satisfies 
$\eta(T)pbp=Tpbp$ for $T\in\cM(D)=\cM(pBp)$ 
and is strictly continuous a map from $\cM(D)$ 
to $\cM(B)$. Since $D=\varphi(A)B\varphi(A)$, 
we can define a \cst-morphism $\lambda$ from 
$A$ into $\cM(D)$ by $\lambda(a)d:=\varphi(a)d$ 
for $a\in A$ and $d\in D$. 
The \cst-morphism $\lambda$ is non-degenerate 
because $A\cdot A = A$ implies that 
$\lambda(A)D=\varphi(A)D=D$ by definition of $D$.
We have seen that non-degenerate $\lambda$ 
extends to a strictly continuous unital \cst-morphism 
$\cM(\lambda)$ from 
$A$ into $\cM(D)$ (in a unique way). 
Define $\cM(\psi)\colon \cM(A)\to \cM(B)$ by 
$\cM(\varphi)(T):=
\eta(\cM(\lambda)(T))\in p\cM(B)p\subseteq \cM(B)$.
We have that $\cM(\varphi)$ is strictly continuous 
because $\eta$ and 
$\cM(\lambda)$ are strictly continuous.
\qed

In particular, 
the non-degeneracy condition 
``$\, \varphi(A)B=B\,$''
is sufficient for the 
existence of $\cM(\varphi)$.

\begin{lem}\label{lem:preserve.corners}
If  $\varphi\colon A\to A$ satisfies 
the above defined property  (CP)
then 
$\varphi_1\colon A\to A_e$
extends to a strictly continuous *-homomorphism
$\cM(\varphi_1)\colon \cM (A)\to  \cM( A_e)$.
\end{lem}
\pf  Let 
$I:=I_\varphi\subseteq  A$ denote the closure
of  $\bigcup_ n (\varphi^n)^{-1}(0)$ 
and suppose that
$I\not=A$, \ie $A_e\not=\{0\}$.
Then $\varphi(I)=I$ and the related class map 
$\psi:=[\varphi]_I$
defines an injective endomorphism of $A/I$.

Notice that epimorphisms $\pi_I\colon A\to A/I$
always extend to 
strictly continuous *-homomorphisms
$\cM(\pi_I)\colon \cM(A)\to \cM (A/I)$.

The *-homomorphism $\varphi\colon A\to A$
extends to a strictly continuous *-homomorphism 
$\cM(\varphi)\colon \cM(A) \to \cM(A)$ if and
only if 
$B:=\varphi(A)+\mathrm{Ann}(\varphi(A))$
contains an approximate unit of 
$A$ by the above mentioned argument.

The natural image $\pi_I(B)$ of
$B=\varphi(A)+\mathrm{Ann}(\varphi(A))$ 
in $A/I$ contains an approximate
unit of $A/I$ and 
$\pi_I(B)\subseteq C + 
\mathrm{Ann}(C)$
for $C:=\pi_I(\varphi(A))= \psi (A/I)$.
It follows that 
$\pi_I \circ \varphi\colon A\to A/I$
extends to strictly continuous
*-homomorphism $\cM(\pi_I\circ \varphi)$
from $\cM(A)$ into $\cM(A/I)$ and that
$\psi\colon A/I \to A/I$ extends
to a strictly continuous *-monomorphism
$\cM(\psi)$ from $\cM(A/I)$ into $\cM(A/I)$.
If we consider any  of the *-homomorphisms
$\varphi_n\colon A\to A_e$ then 
they factorize over 
$A/I$ and there is a natural isomorphism 
$\theta$
from $A_e$
onto $(A/I)_e$, where  
$(A/I)_e$ denotes the inductive limit for
the endomorphism 
$\psi\colon A/I \to A/I$.
It has the following transformations:
$\theta\circ \varphi_n= \psi_n\circ \pi_I$, 
$\psi_{n +1}\circ \psi = \psi_n$
 and
$\psi_n=\sigma \circ \psi_{n+1}$, where 
$\psi_n$ denotes the canonical map
$\psi_{n,\infty}\colon A/I\to A_e$ from the 
$n$'th copy of $A/I$ into the inductive 
$A_e\subseteq A_\infty=\ell_\infty(A)/c_0(A)$.

There is a sequence (constructed below) 
of \cst-subalgebras
$D_1, D_2,\ldots$ of $(A/I)_e$ such that 
$D_n\subseteq  \psi_{n+1}(A/I)$,
$D_n \psi_n(A/I)=\{ 0\}$, and such that
the vector space sum 
$\psi_1(A/I)+ D_1+D_2+\cdots$
contains an approximate unit  of  $(A/I)_e$.

Since $D_nD_m=\{ 0\}$
for
$n\not=m$ and  $D_n \psi_1(A/I)=\{ 0\}$,
this shows  that 
$\psi_1(A/I)+\mathrm{Ann}(\psi_1(A/I), (A/I)_e)$
contains an approximate unit of $(A/I)_e$, and we
can conclude that $\psi_1\colon A/I \to (A/I)_e$
extends to a strictly continuous 
*-monomorphism $\cM(\psi_1)$ from
$\cM( A/I)$ into $\cM((A/I)_e)$.
Since $\theta$ is an isomorphism from $A_e$
onto $(A/I)_e$ and 
$\varphi_1=\theta^{-1}\circ \psi_1 \circ \pi_I$,
we get that 
the the superposition of  strictly continuous
*-homomorphisms
$\cM(\theta)^{-1}\circ \cM(\psi_1)\circ \cM(\pi_I)$
gives a strictly continuous *-homomorphism
from $\cM(A)$ into $\cM(A_e)$ that extends
$\varphi_1\colon A\to A_e$.

The algebras $D_n\subseteq  \psi_{n+1}(A/I)$
can be defined by 
$D_n:= \psi_{n+1}(\mathrm{Ann}(\psi(A/I), A/I))$
using that $\psi_{n+1}\circ \psi =\psi_n$.
Then inductively $\psi_1(A/I)+D_1+\cdots +D_n$
contains an approximate unit of $\psi_{n+1}(A/I)$,
which implies the stipulated existence of an
approximate unit for $(A/I)_e$ in 
$\psi_1(A/I)+D_1+D_2+ \cdots$.
\qed

There are non-equivalent
definitions of crossed products 
by an endomorphism 
in the literature that lead to 
non-isomorphic crossed product \cst-algebras.
Since there are different definitions of 
endomorphism crossed products  
$A\rtimes_\varphi \N$ of $A$ by
the additive semi-group $\N$  of natural numbers,
we describe our definition and notation that 
is inspired by the definitions given by
J.~Cuntz \cite{Cuntz:On}, 
\cite{Cuntz:Symposia38.1982},
W.L.~Paschke
\cite{Paschke80} and P.J.~Stacey 
\cite{Stacey93}.
See \cite{anHuef.Raeburn2012} and 
\cite{Kakariadis2011},
for a general
descriptions of such constructions 
and alternative definitions 
that give different crossed products by $\N$.

The \cst-algebra crossed product 
$A\rtimes_\varphi \N$ associated 
to an endomorphism $\varphi$ of $A$
was  defined by J.\ Cuntz 
in \cite[p.~101]{Cuntz:Symposia38.1982} 
for the special
case where $A$ is unital but 
$\varphi$ is not necessarily
unital. 
It was inspired 
by his special construction in \cite{Cuntz:On} that
showed that $\On$ is a semi-crossed product of 
$M_{n^\infty}$
by the endomorphism  
$\varphi(a):= e\otimes a$ for 
$e:=\mathrm{diag}(1,0,\ldots,0)$.
Since then there where several attempts 
to generalize his construction, 
but not necessarily in a way
that is suitable for our applications.

The generalization of Stacey \cite{Stacey93} 
suffers from his assumption that for \emph{each} 
*-endomorphism $\varphi\colon A\to A$
the natural morphism $\varphi_1\colon A\to A_e$
extends to a *-homomorphism from $\cM(A)$ into
$\cM(A_e)$, or at least to a *-homomorphism
from $\cM(A)$  into $\cM(A_e\rtimes_{\sigma} \Z)$. 
But this is not the case for general
injective *-endomophisms
$\varphi$ of $A$, even if $\varphi$ satisfies the 
above non-degeneracy
condition (ND), \cf{}Example 
\ref{ex:non.extendable.endomorphism}.

With help of the above 
Lemmata 
\ref{lem:corner.image}--\ref{lem:preserve.corners} 
we extend the definition of 
J.~Cuntz's \cite[p.~101]{Cuntz:Symposia38.1982} 
to the non-unital case as follows:

\begin{definition}
Let $\varphi$ be en endomorphism of a 
\cst-algebra $A$ that satisfy
the non-degeneracy property (ND) and 
the corner property (CP).
We define
$A\rtimes_\varphi \N$ to be the 
hereditary \cst-subalgebra
of $A_e\rtimes_\sigma \Z$ that is generated by 
the image $\varphi_1(A)$ of $A$.
\end{definition}
Our
endomorphism $\varphi$ of $A$  
should satisfy
the above discussed non-degeneracy  
property (ND) and 
the corner property (CP).
Indeed, let $B:=A\rtimes _\varphi \N$
denote the crossed 
product induced by an endomorphism
$\varphi$ on $A$ as defined by 
J.\ Cuntz in \cite[p.~101]{Cuntz:Symposia38.1982}.
Knowing that (ND) and (CP) holds one can
\emph{formally} define $B$ as the 
hereditary \cst-subalgebra
of $A_e\rtimes_\sigma \Z$ that is generated by 
the image $\varphi_1(A)$ of $A$.
It is 
even a full hereditary \cst-subalgebra of 
$A_e\rtimes_\sigma \Z$,
because 
$A_e\rtimes_\sigma \Z$ is generated by
$u^n\varphi_1(A)u^m$ for $n,m\in \Z$.

It seems not always to be the case 
that Stacey's version of 
crossed product $A\rtimes_\varphi \N$ 
(see \cite[def.~3.1]{Stacey93} for 
$A\rtimes^1_\varphi \N$)
is naturally isomorphic to the hereditary 
\cst-subalgebra of 
$A_e\rtimes_{\sigma}\Z $ generated by 
$\varphi_1(A)$.
However it is the case when
$\varphi\colon A\to A$
extends to a strictly continuous *-homomorphism
$\cM(\varphi)\colon \cM (A)\to  \cM( A)$ 
-- or equivalently --
that $\varphi_1\colon A\to A_e$ extend to 
$\cM(\varphi_1)\colon \cM (A)\to  \cM( A_e)$,
\cf{}Lemma \ref{lem:preserve.corners}. 
Then, with $p:=\cM(\varphi_1)(1_{\cM(A)})$, 
we have that $p(A_e\rtimes_{\sigma}\Z)p\,$ 
\emph{is the same} as the hereditary 
\cst-subalgebra of $A_e\rtimes_{\sigma}\Z $ 
generated by $\varphi_1(A)$ and is naturally 
isomorphic to Stacey's $A\rtimes _\varphi \N$.

Let us remind the reader of the definition
of a 
$G$-separating action which 
applies to Proposition 
\ref{prop:G.sep.for.endomorph},
and hence also Remark \ref{rem:e.is.dU}
and Lemma 
\ref{lem:equivalent.to.G.separating}.

\begin{definition}[{\cite[def.~5.1]{KirSie1}}]
\label{def:G.separating} 
Suppose that $(A,G,\sigma)$ is a \cst-dynamical 
system with discrete group $G$. The action of 
$G$ on $A$ is \emph{$G$-separating} if
for each $a_1,a_2\in A_+$, $c\in A$ and 
$\varepsilon>0$, there exist $d_1,
d_2\in A$ and $g_1,g_2\in G$ such that 
\begin{equation}\label{InEq.G.separating.fit} 
\| d_j^*a_j d_j-
\sigma_{g_j}(a_j)\|<\varepsilon \quad \text{and} 
\quad \| d_1^*cd_2\| <\varepsilon 
\end{equation}
\end{definition}

\begin{rem}[{\cite[rem.~5.4]{KirSie1}}]
\label{rem:e.is.dU}
Let $(A,G,\sigma)$ a \cst-dynamical system. 
\begin{itemize}
\item[(i)] 
For each 
$a_1,a_2 \in A_+$, $x,d_1,d_2\in A$, 
$g_0,g_1,g_2\in G$ 
and 
$s_1 := d_1U(g_1)$, 
$s_2:=\sigma_{g_0^{-1}}(d_2)U(g_0^{-1}g_2g_2)$,
$c:=x U(g_0)$,
$b_1:=a_1$, and $b_2:=\sigma_{g_0}(a_2)$
the following 
equalities hold:
$$\|  s_j^*a_j s_j - a_j\| = 
\| d_j^*b_j d_j - \sigma_{g_j}(b_j)\|\,
\quad
\textrm{and}\,
\quad
\| s_1^*c s_2\|=\| d_1^*x d_2\|
\,.$$	
\item[(ii)] With $g_0=e$ in (i) the equalities 
reduce to:
$$\|  s_j^*a_j s_j - a_j\| = 
\| d_j^*a_j d_j - \sigma_{g_j}(a_j)\|\,
\quad
\textrm{and}\,
\quad
\| s_1^*c s_2\|=\| d_1^*c d_2\|
\,.$$
\end{itemize}
\end{rem}

\begin{lem}
\label{lem:equivalent.to.G.separating} 
Let $(A,G,\sigma)$ a \cst-dynamical system. 
The following properties are equivalent: 
\begin{itemize} 
\item[(i)]  The action $\sigma$ is $G$-separating, 
\ie for each $a_1,a_2\in A_+$, $c\in A$ and 
$\varepsilon>0$, there exist $d_1,
d_2\in A$ and $g_1,g_2\in G$
satisfying \eqref{InEq.G.separating.fit}.
\item[(ii)] There exists a dense *-subalgebra 
$\cB$ of $A$ that has the properties 
\emph{(1)--(3)}: 
\begin{itemize}
\item[(1)]$\psi(b^*b)\in \cB$ for all 
$\psi\in\Cf_c(0,\infty]_+$ and $b\in \cB$. 
\item[(2)] $\sigma_g(\cB)\subseteq \cB$ for 
all $g\in G$. 
\item[(3)] For each $b_1,b_2\in \cB$, 
$\varepsilon>0$ and $c\in \cB$ there 
exist $d_1,d_2\in A$ 
and $g_1,g_2\in G$ that satisfy the inequalities 
\emph{(\ref{InEq.G.separating.fit})} 
with $a_j:= b_j^*b_j$ ($j=1,2$).
\end{itemize}
\end{itemize}
\end{lem}
\pf Clearly (i)$\Rightarrow$(ii) with $\cB:=A$. 
	
(ii)$\Rightarrow$(i): Define 
$\Co:= \,\{  dU(g) \,;\,\, d\in \cB\,,
\,\, g\in G\, \} \,$ and $\cS:=\Co$.
We use Lemma \ref{lem:cS.control} on the 
above defined $\Co$, $\cS$ and 
$\cF :=\, \{b^*b\,\colon \,\, b\in \cB\}$. 
The property that 
$s_2^*\Co s_1\subseteq  \Co$ for all 
$s_1,s_2\in \cS$ and conditions 
(ii)-(iii) of Lemma \ref{lem:cS.control} 
are trivially satisfied 
using (ii)(1)-(ii)(2). 
Moreover, by (ii)(1), the family $\cF$ is 
invariant under $\varepsilon$-cut-downs, because 
$(b^*b-\varepsilon)_+=\psi(b^*b)^*\psi(b^*b)$ for 
the function $\psi(t):= 
\min\bigl( (t-\varepsilon)_+,\|b\|^2\bigr)^{1/2}
\in \Cf_c(0,\infty]_+$. This implies that also 
condition (i) of Lemma \ref{lem:cS.control} is 
fulfilled by $\cF$, because $\cB$ satisfies 
condition (ii)(3): Take any 
$\varepsilon_0>0$, and $a_j=b_j^*b_j\in \cF$ 
for $j=1,2$. For each $c=xU(g_0)\in \Co$ with 
$x\in \cB, g_0\in G$, and 
$0<\tau \leq \varepsilon\leq \varepsilon_0$ we 
can use (ii)(3) to find elements $
d_1,d_2\in \cB$ and 
$g_1,g_2\in G$ satisfying 
(\ref{InEq.G.separating.fit}) with 
$x,\sigma_{g_0}(a_2),\tau$ in 
place of $c,a_2,\varepsilon$.
Remark \ref{rem:e.is.dU} provides elements 
$s_1,s_2\in \cS$ satisfying \eqref{InEq.general}.
So the pair 
$(a_1,a_2)$ has the matrix diagonalization
with respect to $\cS$ and $\Co $ and 
property
(i) of Lemma 
\ref{lem:cS.control}
holds.

We obtain from Lemma \ref{lem:cS.control} that 
for every $c\in A$, $a_1, a_2\in \cB$, 
$\varepsilon_0/2\geq\varepsilon >0$, and
$\tau>0$, there exist $s_1,s_2\in\cS$ 
satisfying \eqref{InEq.general}.
Using  Remark \ref{rem:e.is.dU} 
we can find (for given $a_j,c,\varepsilon$) 
elements $d_1,d_2\in A$ and $g_1,g_2\in G$
satisfying \eqref{InEq.G.separating.fit}.

\qed

In the following proposition we consider
a dense *-subalgebra   $\cB\subseteq A$  that
is 
\emph{$\varphi$-invariant} -- in the sense  
that $\varphi(\cB)\subseteq \cB$ --
and $\cB$ is a \cst-\emph{local} 
subalgebra  --  in the sense that
$\psi(b^*b)\in \cB$ for 
$b\in \cB$ and $\psi\in \Cf_c(0,\infty]_+$
(see definition in Section \ref{sec2}).
For example,  
$\cB$ can  be an 
algebraic inductive limit 
of an upward directed family of
 \cst-subalgebras of 
 $A$ that is mapped by $\varphi$ into
itself.
\begin{prop}
\label{prop:G.sep.for.endomorph}
Suppose that 
$\varphi$ is an endomorphism of a
\cst-algebra $A$ 
(that is not necessarily injective),
that $\cB\subseteq  A$ is a dense 
*-subalgebra which 
is $\varphi$-invariant,
and that $\cB$ is a \cst-\emph{local} 
subalgebra of $A$.

Let 
$\sigma\colon \Z \to \mathrm{Aut}(A_e)$ 
be the corresponding
action of the integers 
$\Z$ on the inductive limit
$A_e$ of 
the sequence 
$\xymatrix{A \ar[r]^\varphi & A \ar[r]^\varphi & 
A \ar[r]^\varphi & \cdots }\,$.

The following properties $(i)$ and $(ii)$ 
are equivalent: 
\begin{itemize}
\item[(i)] 	For every 
$b_1,b_2, c\in \cB$, and $\varepsilon>0$
there exist 
$k,\,n_1,\,n_2\, \in \N\cup \{0\}$ and elements 
$e_1,\, e_2\,\in A$ such that, for $j\in \{ 1,2\}$,
\begin{equation}\label{m.1.7}
\| \,  e_j^*\varphi^k(b_j^*b_j)\,e_j-
\varphi^{n_j} (b_j^*b_j)\,\|  \,
<\varepsilon\, \quad
\text{and}\,\, \quad 
\| \, e_1^*\varphi^k(c)e_2\,\|  \, 
<\varepsilon\,.
\end{equation}
	
\item[(ii)] The action 
$\sigma\colon \Z\to \mathrm{Aut}(A_e)$ 
of $G:=\Z$  on $A_e$  is  
$G$-separating (Def.\ \ref{def:G.separating}).
\end{itemize}
\end{prop}

\pf
Let $\Co := \bigcup _m \varphi_m(\cB)$.
Since $\cB$ is a dense *-subalgebra of  $A$,
$A_e$ is the closure of
of $\Co $ in $\ell_\infty(A)/c_0(A)$:
$$\Co \subseteq  
\bigcup _m \varphi_m(A)\subseteq  
A_e \subseteq  
\ell_\infty(A)/c_0(A)\,.
$$
Since $\varphi_m(\cB)$ is a 
\cst-local algebra for each
$m\in \N$ and 
$\varphi_m(\cB)\subseteq  \varphi_{m+1} (\cB)$,
the *-subalgebra $\Co $ of $A_e$ is
a dense \cst-local subalgebra of 
$A_e$ that satisfies
$\sigma(\Co )\subseteq  \Co $.

(i)$\Rightarrow$(ii): 
Since 
$\cB$ is a dense *-subalgebra of $A$, 
we may suppose
that the $e_1,e_2\in A$ 
that satisfy the inequalities (\ref{m.1.7})
are actually in $\cB$ itself.

By Lemma \ref{lem:equivalent.to.G.separating}
it suffices to show that, for
$x_1,x_2,y\in \Co \subseteq  A_e$ and 
$\varepsilon>0$
there exists 
$d_1,d_2\in A_e$ and $k_1,k_2\in \Z$  
such that for $j=1,2$,
\begin{equation}
\label{InEq.G.separation.local.subalg.2}
\| d_j^*x_j^*x_jd_j -\sigma^{k_j} (x_j^*x_j) \| 
<\varepsilon
\quad \text{and} \quad 
\| d_1^*yd_2\| <\varepsilon\,.
\end{equation}
Since $\Co $ is the union of the 
family of images 
$\varphi_m(\cB)\subseteq  A_e$
of $\cB$,
there exists $m\in \N$
and $b_1,b_2,c\in \cB$ with 
$\varphi_m(b_j)=x_j$
and $\varphi_m(c)=y$.

We apply the condition in part (i) to 
$(b_1,b_2,c,\varepsilon)$ 
and get $e_1,e_2\in \cB$ and
$k,\,n_1,\,n_2\, \in \N\cup \{0\}$ such that
the inequalities (\ref{m.1.7}) are fulfilled.

Since $\varphi_m$ is a contractive linear map
and  
$\varphi_m\circ \varphi^\ell (a)=  
\sigma^{\ell}(\varphi_m(a))$
for $a\in A$ and $\ell\in \N$,
 we get that 
 $d_j:=\sigma^{-k}( \varphi_m(e_j))$ and 
 $k_j:=n_j-k$
fulfill the  inequalities 
(\ref{InEq.G.separation.local.subalg.2}). 

(ii)$\Rightarrow$(i):
Let $b_1,b_2, c\in \cB$ 
and $\varepsilon>0$.
Since the action of $\Z$ defined by 
$\sigma$ is 
$G$-separating  on $A_e$, 
there exists $d_1,d_2\in A_e$ and 
$k_1,k_2 \in \Z$
such that, for $j=1,2$, 
$$
\| d_j^*\varphi_1(b_j^*b_j)d_j-
 \sigma^{k_j}(\varphi_1(b_j^*b_j)) \| <\varepsilon
\quad \text{and} \quad  
\| d_1^*\varphi_1(c)d_2\| <\varepsilon
\,.$$

Since $\Co $ is a dense *-subalgebra of 
$A_e$ 
we may suppose that 
$d_1,d_2\in \varphi_\ell (\cB)$
for some $\ell\in \N$.
Then there are  $y_1,y_2\in A$ such that
$d_j=\varphi_{\ell} (y_j)$. Since 
$\varphi_{1} = \varphi_{\ell} 
\circ \varphi^{\ell-1}\,$ 
and 
$\sigma^{n}\circ \varphi_{\ell}=
\varphi_\ell\circ \varphi^{n}$
for $n\ge 0$ 
one gets, for $x\in A$,  
$m := \min (0, k_1,k_2)$  and
 $f_j:= \varphi^{-m}(y_j)$  that
\begin{align*}
\varphi_{\ell} 
\bigl(f_i^*\varphi^{\ell-1-m}(x)f_j \bigr) &=
\sigma^{-m}\bigl(d_i^*\varphi_1(x)d_j\bigr),\\
\sigma^{- m}
\bigl( 
\sigma ^{k_j} \circ \varphi_\ell \circ 
\varphi^{\ell-1}(x) 
\bigr)  
&=\varphi_\ell (\varphi^{k_j-m +\ell-1}(x)).
\end{align*}
This gives 
$$ 
\| \varphi_\ell \bigl(\, 
f_j^*\varphi^{\ell-1-m}(b_j^*b_j )f_j - 
\varphi^{k_j-m +\ell-1}(b_j^*b_j)
\,\bigr)\|  < \varepsilon, \ \ 
\| \varphi_{\ell} 
\bigl( f_1^*\varphi^{\ell-1-m}(c)f_2  \bigr) 
\| <\varepsilon
\,.$$
Since 
$\| \varphi_\ell (a) \| =
\lim_{n\to \infty } \| \varphi^n(a)\|\,$
we find sufficiently large $n\in \N$ such that
with 
$$e_j:= \varphi^n(f_j), \quad 
k:= n+ \ell -1-m \quad 
\text{and}\quad   n_j:= n+ k_j-m+\ell-1
$$
the inequalities  (\ref{m.1.7}) are fulfilled.
\qed 

We remind the reader of the notion of properly 
outer actions needed for the next theorem.

\begin{definition}[{\cite{KirSie1}}]
\label{def:properly.outer} 
Suppose that $(A,G,\sigma)$ is a \cst-dynamical 
system and that $G$ is discrete. The action 
$\sigma$ will be called 
\emph{element-wise properly outer}  if, for each 
$g\in G\setminus \{ e\}$, the automorphism 
$\sigma_g$ of $A$ is properly outer in the sense 
of \cite[def.~2.1]{Elliott.prop.outer}, \ie  
$\|\,\sigma_g | I \,-\,\mathrm{Ad}(U)\,\|\,=2$ 
for any  $\sigma_g$-invariant non-zero ideal $I$
 of $A$ and any unitary $U$ in the multiplier 
algebra $\cM(I)$ of $I$. See also 
\cite[thm.~6.6(ii)]{OlePed3}.

An action $\sigma$ is \emph{residually properly 
outer} if, for every $G$-invariant ideal $J\not=A$ 
of $A$, the induced action $[\sigma]_J$ of $G$ on 
$A/J$ is \emph{element-wise} properly outer.
\end{definition}

\begin{thm}\label{thm:endo.cross.spi}
Let  
$\mathcal{B}\subseteq A$, $A_e$, 
$\varphi$ 
and $\sigma$
be as in Proposition 
\ref{prop:G.sep.for.endomorph},
with
endomorphism $\varphi\colon A\to A$
that satisfies the above discussed  
conditions (ND) and (CP).
Suppose that: 
\begin{itemize}
\item[(i)] 	For every $b_1,b_2, c\in \cB$ and 
$\varepsilon>0$ 
there exist 
$k,\,n_1,\,n_2\, \in \N\cup \{0\}$ and elements 
$e_1,\, e_2\,\in A$ such that, for $j\in \{ 1,2\}$,
\begin{equation}\label{m.1.6}
\| \,  e_j^*\varphi^k(b_j^*b_j)\,e_j-
\varphi^{n_j} (b_j^*b_j)\,\|  \, <\varepsilon\,,
\;\;\,\text{and}\,\;\, 
\| \, e_1^*\varphi^k(c)e_2\,\|  \, 
<\varepsilon\,.
\end{equation}
	
\item[(ii)] 	
For every $n\in \N$ and every 
$\sigma$-invariant
closed ideal $J\not=A_e$ of $A_e$ 
the automorphism
$([\sigma]_J)^n$ of $A_e/J$ is
properly outer. 
\end{itemize}
Then $A_e\rtimes_\sigma \Z$ and its
hereditary \cst-subalgebra 
$A\rtimes _\varphi \N$
are strongly purely infinite.
\end{thm}

\pf
It suffices to show that 
$A_e\rtimes_{\sigma}\Z$
is strongly purely infinite, because 
$A\rtimes_\varphi \N$ is naturally isomorphic
to  the (full) hereditary \cst-subalgebra of
$A_e\rtimes_{\sigma} \Z$ that is generated by 
its \cst-subalgebra $\varphi_1(A)$.
Hereditary \cst-subalgebras
of s.p.i.\ algebras are again s.p.i.\ by 
\cite[Prop.~5.11]{KirRorOinf}.
By Proposition \ref{prop:G.sep.for.endomorph},
the condition (i) is equivalent to the 
$G$-separation of the 
action $\sigma$
of $G:=\Z$ on $A_e$ 
generated  by the restriction of the
backward shift on 
$\ell_\infty (A) / c_0(A)$ to $A_e$.
The condition (ii)
says that  the action 
$\sigma$ of $\Z$ on $A_e$
is residually properly outer 
(\cf Definition \ref{def:properly.outer}).
Since every abelian group is exact, the action 
$\sigma$ is exact 
in the sense of \cite[def.~1.2]{Siera2010}.
So all the assumptions  
of \cite[Theorem 1.1]{KirSie1}
are satisfied for $A_e$
and $\sigma
\colon  \Z\to \mathrm{Aut}(A_e)$.
Thus, 
$A_e\rtimes _\sigma \Z$ 
is strongly purely infinite. 
\qed

\textbf{Proof of 
Theorem 
\ref{thm1.4}}:\,
Follows form Theorem 
\ref{thm:endo.cross.spi} as $\cB:=A$ is 
a $\varphi$-invariant \cst-local 
*-subalgebra of $A$.
\qed

\section{Cuntz-Pimsner algebras}

An application of 
Theorem \ref{thm:endo.cross.spi}
to certain special Cuntz-Pimsner 
$\mathcal{O}(\Hi )$ algebras 
is given by the construction below. It is 
implicitly contained in \cite{HarKir}.

Let $C$ be a stable 
$\sigma$-unital \cst-algebra $C$,
and let $h\colon C\to \cM(C)$ be a 
non-degenerate *-homomorphism 
(i.e. $h(C)C=C$) 
that is faithful and satisfies 
$h(C)\cap C=\{ 0\}$. 
Notice that $h$ extends to a 
faithful strictly continuous
unital *-endomorphism $\cM(h)$
of $\cM(C)$.
To simplify notation we denote 
the endomorphism 
$\cM(h)$ of $\cM(C)$ again by $h$, 
unless 
we wish to make an emphasis 
on the difference between 
$\cM(h)$ and $h$.

In the following let 
$\mathcal{H}(h,C)$, or simply 
$\mathcal{H}$,
denote Hilbert bi-module given by
$\mathcal{H}:=C$, with 
right multiplication $b\mapsto ba$,
left-multiplication $b\mapsto h(a)b$  
and Hermitian form 
$\langle c,b\rangle:= c^*b$.
The \cst-algebra $\cL(\mathcal{H})$
of adjoint-able bounded operators on 
$\mathcal{H}$
is here the same as $\cM(C)$.

The closer look in \cite{HarKir}
to the work of Pimsner
\cite{Pimsner97} shows that,
under our special assumptions on 
$h\colon C\to \cM(C)$, the natural 
epimorphism from the
Toeplitz-Pimsner algebra 
$\mathcal{T}(\Hi )$
onto the Cuntz-Pimsner algebra 
$\mathcal{O}(\Hi )$ is an isomorphism,
\emph{and} that $\mathcal{T}(\Hi )$
is isomorphic to a
crossed product 
$A\rtimes _\varphi \N$ in the following
manner:

Consider the algebraic sum
$\cB:= C+h(C)+h^2(C)+\cdots$. 
The algebra $\cB$
is a \cst-local *-subalgebra of $\cM(C)$,
because it is
the algebraic inductive limit of the 
\cst-algebras
$C+h(C)+\ldots+ h^{n}(C)$.
Clearly $h(\cB)\subseteq  \cB$.
Let $A\subseteq  \cM(C)$ 
denote the norm-closure
of $\mathcal{B}$ in $\cM(C)$,
and let $\varphi:= h|A$.
The above stated assumptions on $C$ and 
$h\colon C\to \cM(C)$ 
imply
that the Toeplitz-Pimsner algebra 
$\mathcal{T}(\mathcal{H})$
is naturally isomorphic to the 
semi-group crossed product
$A\rtimes_\varphi \N$.

To prove that 
$\mathcal{T}(\mathcal{H})$
is strongly purely infinite
it suffices to show that 
$\cB\subseteq  A$,
$\varphi$ and 
$(A_e,\Z,k\mapsto \sigma^k)$
satisfy the conditions (i) and (ii)
of 
Theorem \ref{thm:endo.cross.spi}.

It is not possible to prove the
conditions (i) and (ii) with 
the above weak assumptions on 
$h\colon C\to \cM(C)$
that we have introduced so far,
because of an example  
$h\colon C\to \cM(C)$
with $C= \Cf_0(X,\cK)$
(where 
$X:=S^2\times S^2\times \cdots $)
given by M.~R{\o}rdam 
\cite{Rordam.Example}. 
His example has the property that
$\mathcal{O} (\Hi)$ 
is a stable, separable,  simple 
nuclear
\cst-algebra that contains finite 
and infinite projections.
In particular this algebra 
$\mathcal{O} (\Hi)$ 
($\cong \mathcal{T}(\mathcal{H})$)
is not purely infinite.
Therefore we require 
now following stronger properties
(i)--(iv)
for  $h\colon C\to \cM(C)$: 
\begin{itemize}
\item[(i)] $\,h\, $ is a non-degenerate 
*-mono\-mor\-phism.
\item[(ii)] $\, h\, $ is approximately 
unitarily equivalent 
in $\cM (C)$ to its infinite repeat 
  $\delta_\infty\circ h$.
\item[(iii)] Each  
$J\in\mathcal{I} (C)$ is contained
in the closed ideal of 
$C$ generated by $h(J)C$.
\item[(iv)] $\, h\, $ is 
approximately unitarily equivalent
in $\cM (C)$ to $\cM (h)\circ h\,$.
\end{itemize}
Here $\cM (h)\colon \cM(C)\to \cM(C)$ 
denotes the
the strictly continuous extension of 
$h$ to a unital *-mono\-mor\-phism of 
$\cM (C)$ using property (i). 
Recall that the 
\emph{infinite repeat} endomorphism 
$\delta_\infty\colon  \cM(C)\to \cM(C)$
in condition (ii) is unique up to unitary 
equivalence in $\cM(C)$ and
is given by the strictly convergent series
$
\delta_\infty(b):= 
\sum_{n=1}^\infty s_n b s_n^*
$
for $b\in \cM(C)$, where 
$s_1,s_2,\ldots \in \cM(C)$
is a sequence of isometries with 
$\sum_n  s_ns_n^*$  
strictly convergent to $1$.

\begin{definition}
\label{def:residual.weak.Rokhlin}
We say that an action 
$\sigma\colon G\to \mathrm{Aut}(A)$
has the 
\emph{residual weak Rokhlin property}, if the
center $\mathcal{Z}(A^{**})$  of 
$A^{**}$ contains a 
projection 
$P\in  \mathcal{Z}(A^{**})$ that satisfies:
\begin{itemize}
\item[(i)]  $\sigma_g(P)P=0$ for all 
$g\in G\setminus \{ e\}$.
\item[(ii)]
The equation $r(1-q)P=0$ implies $r(1-q)=0$,
if $q,r\in  \mathcal{Z}(A^{**})$ are any 
$\sigma(G)$-invariant 
($A$-)open projections.
\end{itemize}
\end{definition}

Here we extend 
$\sigma_g$ to a normal automorphism
 $\sigma_g\colon A^{**}\to A^{**}$ of $A^{**}$,
and with $q=0$ 
above we obtain the corresponding 
non-residual version of the definition 
which we call the \emph{weak Rokhlin property} 
of the action $\sigma$.

\begin{rem}\label{res.weak.rokhlin}
A comparison of other 
(generalized) Rokhlin properties 
can be found in
\cite[cor.~2.22]{Siera2010}.
Using arguments from the proof of 
\cite[thm.~2.12]{Siera2010} it follows that
topological freeness of \cite[def.~1]{ArchSpiel}
implies topological freeness of 
\cite[def.~1.16]{Siera2010} 
that in turn gives the
weak Rokhlin property 
(Definition \ref{def:residual.weak.Rokhlin})
which again implies the Rokhlin$^*$ 
property \cite[def.~2.1]{Siera2010}
from which we get 
element-wise proper outerness 
(Definition \ref{def:properly.outer}), 
and all properties coincide on 
\emph{commutative} \cst-algebras 
$A\cong \mathrm{C}_0(X)$ and $G$ countable.
It can be easily seen that the proofs of these
results pass to the corresponding versions of 
(generalized) residual Rokhlin properties.
\end{rem}

We outline how above
properties (i)--(iv)
of $h\colon C\to \cM(C)$ 
imply the conditions (i) and (ii)
of Theorem \ref{thm:endo.cross.spi}:
We let  $D:=\varphi_1(C)\subseteq  A_e$, 
then
$$
\varphi_1(\mathcal{B})=
D+\sigma (D) +\sigma^2 (D) +\cdots
$$
and 
$\varphi_1 \circ \varphi= \sigma \circ \varphi_1$
for our above defined automorphism $\sigma$
of $A_e$ associated to $\varphi= h| A$.
It is easy to show (\cf \cite{HarKir}) 
that $A_e$ is the closure of the algebraic sum
$ \sum_{k \in \Z}  \sigma^k (D) $ and that
the closures $J_n$ of  
$\sum _{k\leq n} \sigma^k (D)$
are ideals of  $A_e$ with the property that
$J_n=J_n\sigma^k(D)$ for $k\geq n$, and
$\sigma^n(D)A_e\sigma^{n}(D)=J_n$.
One can use this as a dictionary to translate 
our conditions
on $h\colon C\to \cM(C)$ into conditions on
 $D\subseteq  A_e$ and $\sigma$.
Let $P_0\in  (A_e)^{**}$ denote the 
support projection of the hereditary 
\cst-algebra $DA_eD\subseteq  A_e$. 
Since $DA_eD=J_0$, the projection $P_0$ 
is an open central projection of 
$(A_e)^{**}$.
It is shown in \cite{HarKir} (\,\footnote{\, 
But terminology in \cite{HarKir}
is different, \eg the $\sigma$  used in \cite{HarKir}
is the inverse of our $\sigma$, and 
the $h$ there satisfies weaker assumptions.\,}\,)
that the conditions (i)--(iii) imply that
$P:= P_0- (\sigma^{-1})^{**}(P_0)$
has the properties (i) and (ii) of 
Definition \ref{def:residual.weak.Rokhlin}
for $(A_e,\Z,k\mapsto \sigma^k)$.
Thus conditions 
(i)--(iii) on $h$ imply that the $\Z$-action 
$k \mapsto \sigma^k$
has the residual weak Rokhlin property
of Definition \ref{def:residual.weak.Rokhlin}.
Using Remark \ref{res.weak.rokhlin}
we get property (ii) of Theorem 
\ref{thm:endo.cross.spi}.

It is a fairly elaborate work 
to show that $(\cB, \varphi= h|A )$ satisfy
the inequalities
(\ref{m.1.6}) of Theorem 
\ref{thm:endo.cross.spi}
if $h$ moreover satisfies condition (iv), 
but deep reasonings are not needed.
In this way one can see that (i)-(iv) 
imply the conditions
(i) and (ii) of 
Theorem \ref{thm:endo.cross.spi}.

\medskip

\begin{rems}\label{rems:HHEK.discussion}
Let $C$ be a stable 
$\sigma$-unital \cst-algebra, and 
$h\colon C\hookrightarrow \cM (C)$ 
a *-mono\-mor\-phism.

If $h$ satisfy properties (i)--(iv) and if $C$
 is in addition 
nuclear and separable,
then $\mathcal{O}(\Hi)$
-- build from $h$ --
is a stable separable
nuclear
\cst-algebra that absorbs $\Oinf$ tensorial, 
\ie
\begin{equation}\label{eq:Oinf.absorb}
\mathcal{O}(\Hi ) \cong  
\Oinf \otimes  \mathcal{O}(\Hi )\,.
\end{equation}
We do not know if  the isomorphism 
(\ref{eq:Oinf.absorb})
holds in case that
$h\colon C\to \cM (C)$
satisfies (i)--(iii), but   $C$ is \emph{not} 
nuclear. We  did not find a counter-example for 
the isomorphism  (\ref{eq:Oinf.absorb})
with $h$ satisfying 
only  (i) and (ii).
The property (iii) is  used 
in \cite{HarKir} for the proof of
the residual proper outerness of the 
corresponding $\mathbb{Z}$-action on $A_e$.
The conditions 
(iii), (iv) and the nuclearity
of $C$ play an important role in 
\emph{our} verification of 
the isomorphism (\ref{eq:Oinf.absorb}).
\end{rems}
\begin{rem}
Since many strongly purely infinite nuclear 
\cst-algebras 
are Cuntz-Pimsner algebras of the type  
constructed in \cite{HarKir}
and some of them are stably projection-less, 
our considerations
are also farer going than  for example 
the study of local boundary actions
in \cite{LacaSpi:purelyinf}, 
because reduced crossed products by
local boundary actions 
is very rich of projections
by \cite[lem.~8]{LacaSpi:purelyinf},  but 
there are important amenable 
strongly p.i.\  \cst-algebras
that do not contain any non-zero projection.
\end{rem}

\begin{ex}
\label{ex:non.extendable.endomorphism}
Let $A:= \Cf_0(0,1]$, take 
$g \in A_+$ defined by 
$g(t):= \min(4(t-1/2)_+, t)$
for $t\in [0,1]$. The map
$\varphi\colon\,\, f\in \Cf_0 (0,1] 
\mapsto  f\circ g \in \Cf_0 (0,1]$
is a *-monomorphism of $A$ 
that satisfies condition (ND)
and does not 
extend to
any *-endomorphism $\psi$ of  
$\cM (\Cf_0 (0,1])\cong \Cf_b(0,1]$.
In particular, 
it does not satisfy condition (CD).
\end{ex}

\pf Indeed, $\varphi$ satisfies
property (ND), because $g(t)=t$ for 
$t\in [2/3,1]$.
\emph{Suppose}  that there exists 
is a *-homomorphism 
$\psi\colon \Cf_b(0,1]\to \Cf_b(0,1]$
with $\psi (f)=\varphi(f)=f\circ g$ for
$f\in \Cf_0(0,1]$. Then we would get that
$$f_1(g(t))g(t)=\psi (f_1f_0)(t)=
\psi (f_1)(t) \varphi(f_0)(t)=\psi (f_1)(t)g(t)
$$  
for $t\in (0,1]$, $f_1\in \Cf_b(0,1]$, 
and $f_0(t):=t$.
Thus, $\psi (f_1)(t)=f_1(g(t))$ for all $t > 1/2$.
Since 
$\psi (f_1)$ 
is a bounded continuous function on 
$(0,1]$,  the limit  
$\lim_{\delta \to 0+} \psi (f_1)(1/2+\delta)$ 
exists and is equal to
$\psi (f_1)(1/2)$. It follows that 
$\lim_{\varepsilon \to 0+} f_1(\varepsilon)= 
\lim_{\delta \to 0+} f_1(g(1/2+\delta))$ exists 
for all $f_1\in \Cf_b(0,1]$. The function 
$f_1(t):=\sin (1/t)^2$ 
is in $\Cf_b(0,1]_+$ and has no limit at zero. 
It  contradicts the existence of $\psi $.\qed

\section*{Acknowledgments}
Some of this work was conducted while the 
second named author was at the Fields Institute 
from 2009 to 2012, during a stay of the authors in 
autumn 2011 in Copenhagen 
supported by Prof.~M.~R{\o}rdam 
and 
during the spring of 2014 at the 
Fields Institute during the 
Thematic Program on Abstract Harmonic Analysis, 
Banach and Operator Algebras, 2014 
supported by Prof.~G.~Elliott. 
It is a pleasure for us to forward our thanks to the 
Fields Institute in Toronto and the 
Centre for Symmetry and Deformation at the 
University of Copenhagen. 
In particular, we express our thanks for 
the encouragement and support from 
Prof.~M.~R{\o}dram and Prof.~G.~Elliott. 
We also thank the Mittag-Leffler Institut
for hospitality while this paper
was completed.
The second  authors research was also 
supported by the Australian Research Council.

\providecommand{\bysame}{\leavevmode\hbox
to3em{\hrulefill}\thinspace}
%%%%%%%%%%%%%%%%%%%%%%%%%

\bigskip

\address{Institut f{\"u}r Mathematik,
Humboldt Universit{\"a}t zu Berlin, 
Unter den Linden 6,
D--10099 Berlin, Germany}\\ 
\email{kirchbrg@mathematik.hu-berlin.de}

\address{School of Mathematics \& 
Applied Statistics,
University of Wollongong, 
Faculty of Engineering \& Information Sciences,
2522  Wollongong, Australia}\\
\email{asierako@uow.edu.au}
\end{document}